\documentclass[twoside,12pt]{article}
\usepackage{amssymb}
\usepackage{amsfonts}
\usepackage{amscd}
\usepackage{amsmath}
\usepackage{amsthm}

\advance\topmargin by -\headheight
\advance\topmargin by -\headsep

\newtheorem{theorem}{Theorem}[section]
\newtheorem{proposition}[theorem]{Proposition}
\newtheorem{corollary}[theorem]{Corollary}
\newtheorem{lemma}[theorem]{Lemma}

\theoremstyle{remark}
  \newtheorem{remark}[theorem]{Remark}
  \newtheorem{example}[theorem]{Example}  
\newcommand{\comment}[1]{}

\newcommand \al{\alpha}

\newcommand\ga{\gamma}
\newcommand\de{\delta}

\newcommand\et{\eta}
\renewcommand\th{\theta}
\newcommand\io{\iota}
\newcommand\ka{\kappa}
\newcommand\la{\lambda}
\newcommand\rh{\rho}

\newcommand\si{\sigma}

\newcommand\ph{\varphi}

\newcommand\ps{\psi}
\newcommand\om{\omega}
\newcommand\Ga{\Gamma}

\newcommand\Om{\Omega}

\newcommand\ie{i.e.\ }

\def\RR{\mathbb R}

\def\TT{\mathbb T}
\newcommand\oo{{\infty}}

\renewcommand\o{\circ}
\renewcommand\div{\on{div}}
\newcommand\x{\times}
\newcommand\on{\operatorname}

\newcommand\Ad{\on{Ad}}
\newcommand\ad{\on{ad}}

\newcommand\flux{\on{flux}}

\newcommand\inv{\on{inver}}
\newcommand\ev{\on{ev}}
\newcommand\proj{\on{proj}}

\newcommand\Flux{\on{Flux}}

\newcommand\per{\on{per}}

\newcommand\tr{\on{tr}}

\newcommand\fiber{\on{fiber}}

\newcommand\Quant{\on{Quant}}

\newcommand\Aut{\on{Aut}}
\newcommand\Diff{\on{Diff}}

\newcommand\Pt{\on{Pt}}
\newcommand\Hol{\on{Hol}}

\newcommand\Ker{\on{Ker}}

\newcommand\Ham{\on{Ham}}
\newcommand\Lin{\on{Lin}}

\newcommand\adm{\on{adm}}
\newcommand\SO{\on{SO}}

\newcommand\GL{\on{GL}}
\newcommand\hor{\on{hor}}
\newcommand\eq{\on{eq}}

\newcommand\ex{\on{ex}}

\newcommand\id{\on{id}}

\newcommand\so{\mathfrak s\mathfrak o}
\newcommand\gl{\mathfrak g\mathfrak l}
\newcommand\g{\mathfrak g}

\newcommand\h{\mathfrak h}

\newcommand\X{\mathfrak X}

\newcommand\F{\mathfrak F}

\date{ }
\begin{document}

\title
{Abelian extensions via prequantization}

\author{Cornelia Vizman \\\it \small West University of Timisoara,
Department of Mathematics\\ 
\it\small Bd. V.Parvan 4, 300223-Timisoara, Romania\\
\it \small e-mail: vizman@math.uvt.ro}
\maketitle

\begin{abstract}
We generalize the prequantization central extension of a group of diffeomorphisms
preserving a closed 2--form $\om$, 
to an abelian extension of a group of diffeomorphisms preserving
a closed vector valued 2--form $\om$ up to a linear isomorphism
($\om$--equivariant diffeomorphisms).
Every abelian extension of a simply connected Lie group can be obtained 
as the pull-back of such a prequantization abelian extension.
\end{abstract}



\section{Introduction}\label{intro}

Not every infinite dimensional Lie algebra can be integrated to a Lie group.
Moreover, given a Lie group $G$ with Lie algebra $\g$, not every abelian 
Lie algebra extension of $\g$ can be integrated to an abelian Lie group extension of $G$.
The two obstructions to the integrability of such abelian extensions are  described
in terms of a flux homomorphism involving $\pi_1(G)$ and a period homomorphism 
involving $\pi_2(G)$ \cite{N04} (see the appendix).

Taking this into account, geometric constructions of Lie group extensions in infinite dimensions
are important. A prototype is the prequantization extension \cite{K70}\cite{S70} associated
to a prequantizable symplectic manifold $(M,\om)$: 
\begin{equation*}
1\to S^1\to \Quant(P,\th)_0\to\Ham(M,\om)\to 1.
\end{equation*}
Here $P$ is a principal circle bundle over $M$ with connection $\th$ and curvature $\om$.
The identity component of the group of quantomorphisms, \ie connection preserving automorphisms,
is a 1-dimensional central extension of the group of hamiltonian diffeomorphisms of $M$.
This is a Lie group extension when $M$ is compact \cite{RS81}.
Central Lie group extensions of a Lie group $G$ can be obtained as pull-back 
of the prequantization extension, given a hamiltonian $G$-action on the prequantizable symplectic manifold $(M,\om)$.

This construction works more generally for closed vector valued 2-forms $\om$ on 
a smoothly paracompact manifold $M$, possibly infinite dimensional, with values in a Mackey complete
locally convex vector space $V$, and with discrete period group $\Ga\subset V$,
\ie for prequantizable forms. 
The identity component of the  "quantomorphism group" $\Diff^{\eq}(P,\th)^A$ of the principal bundle $q:P\to M$ with abelian structure group $A=V/\Ga$ and principal $A$-action $\rho$ is an abelian extension of the "hamiltonian group" $\Diff_{\ex}(M,\om)$
\begin{equation*}
1\to A\to \Diff(P,\th)^A_0\to\Diff_{\ex}(M,\om)\to 1.
\end{equation*}
 
In the infinite dimensional setting, the prequantization extension is not a Lie group extension,
but the pull-back provides Lie group extensions of Lie groups 
with hamiltonian actions $(M,\om)$ \cite{NV03}.
For instance, all the central extensions of the group of exact volume preserving diffeomorphisms
can be obtained from its hamiltonian action on the  non-linear Grassmannian 
of codimension two submanifolds of $M$ \cite{I96} \cite{HV04}.

Our main goal is to construct a prequantization abelian extension associated to a prequantizable 
form $\om\in\Om^2(M,V)$:
\begin{equation*}
1\to A\to \Diff^{\eq}(P,\th)^A_0\to\Diff_{\ex}^{\eq}(M,\om)\to 1.
\end{equation*}
The identity component of the group $\Diff^{\eq}(P,\th)^A$ of equi-quantomorphisms is an abelian extension
of the group $\Diff_{\ex}^{\eq}(M,\om)$ of equi-hamiltonian diffeomorphisms.
We show that, in the context of abelian extensions,
it plays the same role as the prequantization extension in the context of central extensions.
Moreover, Theorem \ref{th1} shows that every abelian Lie group extension of a simply connected Lie group
can be obtained as pull-back of such a prequantization abelian extension.

An $\om$-equivariant $G$-action involves an action on $M$ and one on $V$.
The group 
\[
\Diff^{\eq}(M,\om):=\{(\ph, u)\in\Diff(M)\x\GL(V):\ph^*\om= u\cdot\om\}
\]
of $\om$-equivariant diffeomorphisms replaces the group $\Diff(M,\om)$ 
of $\om$-invariant diffeomorphisms.
The flux homomorphism integrating the infinitesimal flux homomorphism
\[
\flux:\X(M,\om)\to H^1(M,V),\quad\flux(\et)=[i_\et\om] 
\]
to $\Diff(M,\om)$ can be extended to a flux 1-cocycle on $\Diff^{\eq}(M,\om)$,
integrating the infinitesimal flux 1-cocycle 
\[
\flux^{\eq}:(\et,\ga)\in\X^{\eq}(M,\om)\mapsto [q^*i_\et\om-\ga\cdot\th]\in
H^1(P,V).
\]

The quantomorphisms are automorphisms of the principal bundle $P$, so the (gauge)  abelian extension 
\begin{equation*}
1\to C^\oo(M,A)\to \Diff(P)^A\to\Diff(M)_{[P]}\to 1,
\end{equation*}
contains the prequantization extension.
In contrast to the quantomorphisms, 
the equi-quantomorphisms
are not $A$-equivariant. 
Consequently,  the gauge extension has to be enlarged. 
An extension
of the group $\Diff(M)_{[P]}$ of diffeomorphisms preserving the isomorphism class of $P$
by the non-abelian group of almost $A$-invariant maps 
\[
C^\oo_{A}(P,A)
=\{f\in C^\oo(P,A): \forall a\in A, f^{-1}(f\o\rh(a)) \text{ constant on $P$}\}
\]
contains the prequantization abelian extension.

The plan of the paper is the following:  Sect.~\ref{s2} collects known results concerning central extensions;
Sect.~\ref{s3} introduces the flux 1-cocycle, a generalization of the symplectic
flux homomorphism; Sect.~\ref{s4} contains the construction of the prequantization
abelian extension; in Sect.~\ref{s5} we provide abelian Lie group extensions
via prequantization; the case of an exact 2-form $\om$, when the extensions are given by group cocycles, 
is treated in Sect.~\ref{s6}; Sect.~\ref{s7} is devoted to examples.


\section{Central extensions}\label{s2}

The prequantization extension associated
to a prequantizable symplectic manifold $(M,\om)$ is  
\begin{equation}\label{zero}
1\to S^1\to \Quant(P,\th)_0\to\Ham(M,\om)\to 1.,
\end{equation}
for $P$ a principal circle bundle over $M$ with connection $\th$ and curvature $\om$ \cite{K70} \cite{S70}.
In this section we present 
the prequantization extension in the more general case of a vector valued closed 2--form, following \cite{NV03}. 

We consider a connected manifold $M$, possibly infinite dimensional,
and we assume it is smoothly paracompact, in order to classify principal bundles over $M$ \cite{B}. 
Let $\om$ a closed 2--form on $M$ with values
in the Mackey complete locally convex space $V$.
The period group $\Ga$ of $\om$ is the image of the homomorphism
$[\si]\in H_2(M,\RR)\mapsto\int_{[\si]}\om\in V$.
We assume that the period group $\Ga\subset V$ is discrete, \ie $\om$ is prequantizable. 
Let $\Diff(M,\om)$ be the group 
of $\om$--invariant diffeomorphisms and $\X(M,\om)$ the Lie algebra of $\om$--invariant vector fields. 
The {\it infinitesimal flux homomorphism}
\begin{equation}\label{flu}
\flux:\X(M,\om)\to H^1(M,V),\quad\flux(\et)=[i_\et\om], 
\end{equation}
can be integrated to the {\it flux homomorphism} \cite{NV03}
\begin{gather}\label{elev}
\Flux:\Diff(M,\om)_0\to H^1(M, V)/H^1(M,\Ga)\nonumber\\
\Flux(\ph)=\int_0^1[i_{\de^r\ph_t}\om] dt\mod H^1(M,\Ga)
=\int_0^1[i_{\de^l\ph_t}\om] dt\mod H^1(M,\Ga),
\end{gather}
where $\ph_t$ is any smooth curve  in $\Diff(M,\om)$ connecting the identity and $\ph$,
and $\de^l$, $\de^r$ denote the left, right logarithmic derivative (defined in the Appendix).
The identity component $\Diff_{\ex}(M,\om)$ of the kernel of the flux homomorphism Flux is called the group
of {\it exact $\om$-invariant diffeomorphisms}. 
As in \cite{NV03}, by the identity component $G_0$ of a group $G$ of diffeomorphisms of $M$
we understand the normal subgroup of those diffeomorphisms of $G$ connected to the identity by
a curve $\ph:I\to G$ such that the map $(t,x)\mapsto (\ph(t)x,\ph(t)^{-1}x)$ is smooth.

For a symplectic manifold $(M,\om)$, the {\it symplectic flux homomorphism}
$S_\om$ involves a subgroup, 
$\Pi\subseteq H^1(M,\Ga)$, called the flux subgroup, known to be discrete \cite{On06},
so 
$$S_\om:\Diff(M,\om)_0\to H^1(M,\RR)/\Pi,\quad S_\om(\ph)=\int_0^1[i_{\de^r\ph_t}\om] dt\mod\Pi.$$
In this case the group of hamiltonian diffeomorphisms 
$\Ham(M,\om)$, defined as the group of those diffeomorphisms of $M$ which are endpoints of hamiltonian isotopies,
 is the kernel of the symplectic flux $S_\om$ \cite{B78}, 
and it coincides with $\Diff_{\ex}(M,\om)$ \cite{NV03}.
By abuse of language we call $\Diff_{\ex}(M,\om)$ the "hamiltonian" group
even if $\om$ is not symplectic.


\paragraph{Extensions of Lie algebras of vector fields.}
In this paragraph we present several Lie algebras of vector fields on a principal bundle: 
projectable vector fields,
$A$-invariant vector fields, and infinitesimal quantomorphisms,
as extensions of Lie algebras of vector fields on the base.

Let $A$ be the abelian Lie group $V/\Ga$, where $\Ga$ is the discrete period group of 
the closed vector valued 2-form $\om$.
Since $M$ is smoothly paracompact,
there exists a principal $A$-bundle $q:P\to M$ and a principal connection
1--form $\th\in\Om^1(P,V)$ with curvature $\om\in\Om^2(M,V)$. 
The principal $A$--action is denoted by  $\rh$, and the infinitesimal action by $\dot\rh:V\to\X(P)$.
In particular $d\th=q^*\om$ and $i_{\dot\rh(v)}\th=v$ for all $v\in V$. 

To every function $h\in C^\oo(P, V)$ one associates the vertical 
vector field $\dot\rh(h)$ on $P$ by
$\dot\rh(h)(y)=\dot\rh(h(y))(y)$. 
We endow $C^\oo(P, V)$ with a Lie
bracket such that the injective mapping 
$
\dot\rh: C^\oo(P, V)\to\X(P)
$
becomes a
Lie algebra homomorphism for the opposite Lie bracket on $\X(P)$. 
This leads to 
\begin{equation}\label{[]}
[h_1,h_2]=L_{\dot\rh(h_2)}h_1
-L_{\dot\rh(h_1)}h_2,
\end{equation} 
because 
$i_{[\dot\rh(h_1),\dot\rh(h_2)]}\th
=L_{\dot\rh(h_1)}h_2-L_{\dot\rh(h_2)}h_1$.

A vector field $\xi\in\X(P)$ is called {\it projectable} if it is $q$-related to a vector field $\et\in\X(M)$, and we denote $\et=q_*\xi$. Projectable vector fields can be characterized by
$Tq\o \xi\o\rh(a)=Tq\o\xi$ for all $a\in A$. 
Every vertical vector field is projectable and we get an exact sequence of Lie algebras
\begin{equation}\label{proj}
0\to C^\oo(P, V)\stackrel{\dot\rh}{\to}\X_{\proj}(P)\stackrel{q_*}{\to}\X(M)\to 0.
\end{equation} 

The pull-back by $q$ maps $C^\oo(M,V)$ into an abelian Lie subalgebra of $C^\oo(P,V)$ because the bracket (\ref{[]}) 
on pull-back functions $q^*f$ for $f\in C^\oo(M,V)$ vanishes.
The Lie algebra $\X(P)^A$ of {\it $A$--invariant vector fields} on $P$
(infinitesimal automorphisms of $P$) consists of vector fields $\xi$ such that $\rh(a)^*\xi=\xi$ for all $a\in A$,
or equivalently $L_{\dot\rh(v)}\xi=0$ for all $v\in V$. 
Restricting (\ref{proj}) to the Lie algebra $\X(P)^A$, we obtain an abelian Lie algebra extension 
\begin{equation}\label{ainv}
0 \to C^\oo(M, V)\stackrel{\dot\rh}{\to}\X(P)^A\stackrel{q_*}{\to}\X(M)\to 0.
\end{equation}
Its cohomology class is given by the curvature form $\om$ on $M$ viewed as a 
Lie algebra 2--cocycle on $\X(M)$ with values in the $\X(M)$--module $C^\oo(M,V)$. 

An  {\it infinitesimal quantomorphism} is an infinitesimal connection preserving automorphism of $P$.
The Lie algebra of infinitesimal quantomorphisms can be expressed as
$$
\X(P,\th)^A
=\{\xi\in\X_{\proj}(P):L_\xi\th=0\}
$$ 
because if $q_*\xi=\et$ and $L_\xi\th=0$, then $[\xi,\dot\rh(v)]=0$ for all $v\in V$.
Indeed, $[\xi,\dot\rh(v)]$ is a vertical vector field ($q$-related to $[\et,0]=0$) and $i_{[\xi,\dot\rh(v)]}\th=L_\xi v=0$, so $\xi\in\X(P)^A$.
Restricting (\ref{ainv}) further to 
$\X(P,\th)^A$, we get the prequantization central Lie algebra extension 
\begin{equation}\label{quan}
0\to V\stackrel{\dot\rh}{\to} \X(P,\th)^A\stackrel{q_*}{\to}\X_{\ex}(M,\om)\to 0.
\end{equation}
Indeed, let $\xi\in\X(P,\th)^A$. 
Both $\xi$ and $\th\in\Om^1(P,V)$ being $A$--invariant,
the function $i_\xi\th\in C^\oo(P,V)$ is $A$--invariant too,
hence it descends to a function $q_*i_\xi\th$ on $M$.
Now $L_\xi\th=0$ and $q_*\xi=\et$ imply $i_\et\om=d(-q_*i_\xi\th)$, so $\xi$ is $q$-related to the hamiltonian vector field $\et$.
On the other hand $\dot\rh(h)\in\X(P,\th)^A$ implies $0=L_{\dot\rh(h)}\th=dh$, so all vertical infinitesimal quantomorphisms are of the form $\dot\rh(v)$, $v\in V$.
The cohomology class  
describing this extension is the class of the $V$-valued Lie algebra 
2--cocycle on $\X_{\ex}(M,\om)$ given by $(\et_1,\et_2)\mapsto -\om(\et_1,\et_2)(x_0)$,
$x_0\in M$.

\paragraph{Geometric prequantization.}
A function $f\in C^\oo(M,V)$ is a hamiltonian function for the  vector field $\et_f$ if $i_{\et_f}\om=df$, so  
the hamiltonian functions on $M$ have to be
constant along the leaves of $\Ker\om\subset TM$. They form the
subspace of {\it admissible functions} $C_{\adm}^\oo(M,V)$.
A hamiltonian
vector field associated to such an admissible function $f$ can be
determined only up to a section in $\Ga(\Ker\om)$. 

\begin{remark}\label{corr}
The linear map
\begin{equation}\label{cosp}
\xi\in\X(P,\th)^A\mapsto -q_*i_\xi\th\in C_{\adm}^\oo(M,V)
\end{equation}
is surjective with kernel $\Ga(\Ker\om)^{\hor}$. 
By definition, the horizontal lift of a vector field 
$\et\in\X(M)$ is the unique vector
field $\et^{\hor}$, $q$-related to $\et$, satisfying $i_{\et^{\hor}}\th=0$.
Given an admissible function $f$, there exists a hamiltonian vector field $\et_f$, and the vector field $\xi_f=\et_f^{\hor}-\dot\rh(q^*f)$ is an infinitesimal quantomorphism with $i_{\xi_f}\th=-q^*f$. 

In the symplectic case $\Ker\om=0$, so we have
$C_{\adm}^\oo(M)=C^\oo(M)$ and the hamiltonian vector field $\et_f$ is
uniquely determined by its hamiltonian function $f$.
The linear map (\ref{cosp}) is a bijection
with inverse 
\begin{equation}\label{inv}
f\in C^\oo(M)\mapsto \xi_f:=\et_f^{\hor}-(q^*f)E\in\X(P,\th)^A,
\end{equation}
$E=\dot\rh(1)$ denoting the infinitesimal generator 
of the circle action on $P$.
This is the symplectic prequantization, in the construction due to Souriau
\cite{S70}.
\end{remark}


\paragraph{Extensions of diffeomorphism groups.}
There is a bijection between $C^\oo(P,A)$ and the space $C^\oo_{\fiber}(P,P)$ of fiber preserving smooth maps
which associates to $f\in C^\oo(P,A)$  the map
$$
\rh(f):y\in P\mapsto\rh(f)(y)=\rh(y,f(y))\in P.
$$ 
The composition on $C^\oo_{\fiber}(P,P)$ determines a monoid structure on $C^\oo(P,A)$, namely
\begin{equation}\label{dot}
(f_1\cdot f_2)(y)=f_1(\rh(y,f_2(y)))f_2(y). 
\end{equation}
The image by $\rh$ of the group of invertible elements $C^\oo(P,A)_{\inv}$
in $C^\oo(P,A)$
is the group of fiber preserving diffeomorphisms of $P$.

The group $\Diff_{\proj}(P)$ of {\it projectable diffeomorphisms} is the group of diffeomorphisms of $P$ which map
fibers to fibers, i.e. those $\ps\in\Diff(P)$
such that $q\o\ps=\ph\o q$ for some $\ph\in\Diff(M)$. 
Projectable diffeomorphisms of $P$ can be characterized by $q\o\ps\o\rh(a)=q\o\ps$ for all $a\in A$.
We write $\ph=q_*\ps$ and the diffeomorphism $\ph$ belongs to $\Diff(M)_{[P]}$,
the group of diffeomorphisms
preserving the isomorphism class $[P]$ of the principal bundle $P$.
The exact sequence of groups
\begin{equation}\label{star}
1\to C^\oo(P,A)_{\inv}
\stackrel{\rh}{\to}\Diff_{\proj}(P)
\stackrel{q_*}{\to}\Diff(M)_{[P]}\to 1
\end{equation}
is the global version of (\ref{proj}).

Let $\Diff(P)^A$ be the group of $A$--equivariant diffeomorphisms of $P$, \ie 
the group of automorphisms of the principal bundle $P$.
An abelian extension (the gauge extension) is obtained by restricting \eqref{star} to the subgroup $\Diff(P)^A\subset\Diff_{\proj}(P)$:
\begin{equation}\label{a}
1\to C^\oo(M,A)\stackrel{\rh}{\to}\Diff(P)^A
\stackrel{q_*}{\to}\Diff(M)_{[P]}\to 1,
\end{equation}
with infinitesimal version the abelian Lie algebra extension (\ref{ainv}). 

A {\it quantomorphisms} is a connection preserving automorphism of $P$.
The group of quantomorphisms can be expressed also as
$$
\Diff(P,\th)^A
=\{\ps\in\Diff_{\proj}(P):\ps^*\th=\th\}.
$$
The quantomorphisms of $P$ descend to holonomy preserving diffeomorphisms on $M$. 
Denoting by $h(\ell)\in A$ the holonomy around a loop $\ell$ in $M$ for the principal connection $\th$, let
\begin{equation*}
\Hol(M,\om)=\{\ph\in\Diff(M):\forall\ell\in C^\oo(S^1,M),h(\ph\o\ell)=h(\ell)\}
\end{equation*}
be the group of {\it holonomy preserving diffeomorphisms}. It is a subgroup of the group $\Diff(M,\om)$ of $\om$-preserving diffeomorphisms.

A central extension can be obtained by a 
further restriction of \eqref{a} to the group of quantomorphisms \cite{NV03}: 
\begin{equation}\label{1111}
1\to A\stackrel{\rh}{\to}\Diff(P,\th)^A\stackrel{q_*}{\to}\Hol(M,\om)\to 1,
\end{equation}
An argument can be given using Proposition \ref{appendix} in the Appendix:
for $f\in C^\oo(P,A)$, $\rh(f)^*\th=\th$ if and only if $\de^l(f)=0$, so $f$ is a constant $\in A$.
Passing to connected components of the identity we get the prequantization central extension
\begin{equation}\label{2222}
1\to A\stackrel{\rh}{\to}\Diff(P,\th)^A_0\stackrel{q_*}{\to}\Diff_{\ex}(M,\om)\to 1,
\end{equation}
the analogue of \eqref{zero} 
in the symplectic setting.
 

\paragraph{Geometric construction of central extensions.}

Let $G$ be a connected Lie group, $\la$ a smooth $G$-action on $M$ with infinitesimal action $\dot\la$, and
$\om\in\Om^2(M,V)$ a $G$--invariant closed 2--form with discrete period group $\Ga\subset V$. 
The $G$--action $\la$ 
is called a {\it hamiltonian action} if $i_{\dot \la(X)}\om\in\Om^1(M,V)$ is exact for all $X\in\g$,
so it gives a group homomorphism $\la:G\to\Diff_{\ex}(M,\om)$. 

In Theorem 3.4 from \cite{NV03} is proven that, given 
a hamiltonian action $\la$ of a connected Lie group $G$ on $(M,\om)$,
there exists a central Lie group extension $\hat G$ of $G$ by $A=V/\Ga$ 
and a smooth $\hat G$--action on $(P,\th)$ by 
quantomorphisms, lifting the $G$--action. 
The central extension of $G$ is a pull-back of the prequantization central extension (\ref{2222}) and a corresponding Lie algebra cocycle on $\g$ is $(X,Y)\mapsto-\om(\dot \la(X),\dot \la(Y))(x_0)$, where $x_0\in M$ is fixed.
The manifold structure on $\hat G$ is obtained from the pull-back bundle of $P$ by an orbit map of $G$ on $M$.


\section{Flux 1--cocycle}\label{s3}

Let $M$ be a connected smoothly paracompact manifold, possibly infinite dimensional, 
$V$ a Mackey complete locally convex space, and
$\om\in\Om^2(M,V)$ prequantizable, \ie closed 2--form with discrete period group.
In this section we introduce the group $\Diff^{\eq}(M,\om)$ of $\om$-equivariant diffeomorphisms
as an isotropy subgroup of $\Diff(M)\x\GL(V)$,
and we define a flux 1-cocycle, extending the flux homomorphism \eqref{elev} on $\Diff(M,\om)$.
Properties of 1-cocycles can be found in the appendix.
 
The Lie algebra of {\it $\om$--equivariant vector fields} on $M$,
\begin{equation*}
\X^{\eq}(M,\om):=\{(\et,\ga)\in\X(M)\x\gl(V):L_\et\om=\ga\cdot\om\},
\end{equation*}
is the stabilizer of $\om$ under the representation $(\et,\ga)\cdot\om=-L_\et\om+\ga\cdot\om$
of the direct product Lie algebra $\X(M)\x\gl(V)$.
We take the negative sign convention on Lie algebras of vector fields
(so $\X(M)$ is the Lie algebra of the group of diffeomorphisms of $M$), so
\begin{equation}\label{lb}
[(\et_1,\ga_1),(\et_2,\ga_2)]=(-[\et_1,\et_2],\ga_1\ga_2
-\ga_2\ga_1).
\end{equation}
The Lie algebra of $\om$--invariant vector fields 
is a Lie subalgebra of $\X^{\eq}(M,\om)$ for the inclusion 
\begin{equation*}
\io:\X(M,\om)\to\X^{\eq}(M,\om),\quad\io(\et)=(\et,0). 
\end{equation*}

\begin{remark}\label{vzero}
Let $V_0$ be the closure of
the image of $\om:TM\x_M TM\to V$.
Given an $\om$-equivariant vector field $(\et,\ga)$, the restriction of $\ga$ to $V_0$ is uniquely determined by $\et$.
When $V_0=V$, then $\ga$ is determined by $\et$ and one can identify the Lie algebra $\X^{\eq}(M,\om)$ with its projection on the first factor: $\{\et\in\X(M):\exists\ga\in\gl(V)\text{ s.t. }L_\et\om=\ga\cdot\om\}$.
When $V_0\ne V$, an $\om$-invariant vector field $\et$ can determine other $\om$-equivariant vector fields beside $(\et,0)$, namely $(\et,\ga)$ with $\ga|_{V_0}=0$.

The period group $\Ga$ is a subgroup of $V_0$ and,
for any $(\et,\ga)\in\X^{\eq}(M,\om)$, the restriction of $\ga$ to $\Ga$
is trivial. Indeed, $\ga(\int_\si\om)=\int_\si L_\et\om=0$ for any 2--cycle $\si$ in $M$.
In particular $\ga=0$ if $V$ is generated by $\Ga$, so $\io(\X(M,\om))=\X^{\eq}(M,\om)$ in this case.
For a closed $\RR$--valued 2--form $\om$, the Lie algebra of $\om$--equivariant 
vector fields is strictly bigger than the Lie algebra of $\om$--invariant
vector fields if and only if $\Ga=0$, \ie $\om$ is exact. 
\end{remark}

\paragraph{Infinitesimal flux cocycle.}
The projection on the second factor, 
\begin{equation}\label{acti}
(\et,\ga)\in\X^{\eq}(M,\om)\mapsto\ga\in\gl(V),
\end{equation}
is a Lie algebra homomorphism, so $V$ becomes a $\X^{\eq}(M,\om)$--module in a natural way. 
The action of the image of $\io$ is trivial.

Since the period group $\Ga\subset V$ of $\om$ is discrete, 
we consider again a principal $A=V/\Ga$--bundle $q:P\to M$ with connection
form $\th\in\Om^1(P,V)$ and curvature form $\om\in\Om^2(M,V)$. 

\begin{proposition}\label{si} 
The linear map 
\begin{equation}\label{seq}
\flux^{\eq}:(\et,\ga)\in\X^{\eq}(M,\om)\mapsto [q^*i_\et\om-\ga\cdot\th]\in
H^1(P,V)
\end{equation}
is a Lie algebra 1--cocycle for the natural
$\X^{\eq}(M,\om) $--module structure on $H^1(P,V)$ induced
by the $\X^{\eq}(M,\om)$--action \eqref{acti} on $V$. 

Its cohomology class $[\flux^{\eq}]\in H^1(\X^{\eq}(M,\om),H^1(P,V))$ is independent of the choice of the connection $\th$. 
\end{proposition}

\begin{proof}
The 1--form $q^*i_\et\om-\ga\cdot\th$ on $P$ is closed for
any $\om$--equivariant vector field $(\et,\ga)$ because
$L_\et\om=\ga\cdot\om$ and $q^*\om=d\th$.
For $(\et_1,\ga_1),(\et_2,\ga_2)\in\X^{\eq}(M,\om)$,
\begin{align*}
\flux^{\eq}([(\et_1,\ga_1),&(\et_2,\ga_2)])
=-[q^*i_{[\et_1,\et_2]}\om+(\ga_1\ga_2-\ga_2\ga_1)\cdot\th]\\
&=[d(q^*\om(\et_1,\et_2))+q^*i_{\et_2}L_{\et_1}\om-q^*i_{\et_1}L_{\et_2}\om
-\ga_1\ga_2\cdot\th+\ga_2\ga_1\cdot\th]\\
&=[q^*i_{\et_2}(\ga_1\cdot\om)-\ga_1\ga_2\cdot\th]-[q^*i_{\et_1}(\ga_2\cdot\om)-\ga_2 \ga_1\cdot\th]\\
&=\ga_1\cdot\flux^{\eq}(\et_2,\ga_2)-\ga_2\cdot\flux^{\eq}(\et_1,\ga_1),
\end{align*}
so the 1--cocycle condition for $\flux^{\eq}$ is satisfied for
the natural $\X^{\eq}(M,\om)$--action.

Two connection 1--forms on $P$ differ by the pull-back $q^*\al$ of a closed $V$--valued 1--form $\al$ on $M$. Then the corresponding flux 1--cocycles differ by the linear map $(\et,\ga)\mapsto \ga\cdot[q^*\al]$, which is a 1--coboundary on the Lie algebra $\X^{\eq}(M,\om)$. 
\end{proof}

The 1--cocycle $\flux^{\eq}$ is called {\it the infinitesimal
flux 1--cocycle}.
Its kernel, denoted by $\X_{\ex}^{\eq}(M,\om)$, is a Lie
subalgebra of $\X^{\eq}(M,\om)$ and is called the Lie algebra of 
{\it equi-hamiltonian
vector fields} on $(M,\om)$. 
We say that $h\in C^\oo(P,V)$ is an {\it equi-hamiltonian function} 
for the equi-hamilto\-nian vector field $(\et,\ga)\in\X^{\eq}_{\ex}(M,\om)$ if 
\begin{equation}\label{defham}
q^*i_\et\om-\ga\cdot\th=dh.
\end{equation}
The equi-hamiltonian function doesn't determine uniquely the equi-hamiltonian vector field,
unless $\ker\om=0$.

\begin{remark} 
The infinitesimal flux
homomorphism \eqref{flu} and the infinitesimal flux 1-cocycle \eqref{seq} are related by
$\flux^{\eq}\o\io=q^*\o\flux$, hence the inclusion $\io$ descends to an inclusion 
$\io:\X_{\ex}(M,\om)\to\X_{\ex}^{\eq}(M,\om)$
of the ideal of exact $\om$--invariant vector fields (hamiltonian vector fields
when $\om$ is symplectic) in the Lie algebra of equi-hamiltonian vector fields.
\end{remark}

\paragraph{Admissible functions.}
Not every smooth $V$--valued function on $P$ can play the role
of an equi-hamiltonian function. We denote by $C^\oo_{\adm}(P,V)$
the space of all possible equi-hamiltonian functions, also called {\it admissible functions}.

\begin{proposition}\label{adm}
If $h\in C_{\adm}^\oo(P,V)$ is an equi-hamiltonian function
for the equi-hamil\-tonian vector field $(\et,\ga)$,
then $L_{\dot\rh(v)}h=-\ga(v)$ for all  $v\in V$.
There exists a group homomorphism $\bar\ga:A\to V$ satisfying
$\bar\ga\o\exp=\ga$, with $\exp:V\to A$ the canonical projection,
such that for all $a\in A$ holds
$h-h\o\rh(a)=\bar\ga(a)$.
\end{proposition}

\begin{proof}
Since $h$ is an equi-hamiltonian function for the equi-hamiltonian 
vector field $(\et,\ga)$, we have that 
$dh=q^*i_\et\om-\ga\cdot\th$.
Then from $L_{\dot\rh(v)}h=i_{\dot\rh(v)}dh=-\ga(v)$ we get the first identity. 

From $q\o\rh(a)=q$ and $\rh(a)^*\th=\th$ it follows that $\rh(a)^*dh=dh$,
so $h-h\o\rh(a)$ is a constant function on the connected manifold  $P$.
This ensures the existence of a group homomorphism $\bar\ga:A\to V$ satisfying the identity
$h-h\o\rh(a)=\bar\ga(a)$ for all $a\in A$.
From that we easily get that $\bar\ga\o\exp\in\gl(V)$.
To show that $\bar\ga\o\exp=\ga$, we differentiate at $t=0$ the identity 
$h-h\o\rh(\exp tv)=\bar\ga(\exp tv)$ for $v\in V$,
and we obtain that
$L_{\dot\rh(v)}h=-T_1\bar\ga(v)$. This gives $\ga=T_1\bar\ga=T_0(\bar\ga\o\exp)=\bar\ga\o\exp$.
\end{proof}

The space of {\it almost $A$--invariant functions} is
\begin{align}\label{capv}
C^\oo_A(P,V)
=\{h\in C^\oo(P,V):\forall a\in A,h-h\o\rh(a)=\text{ constant on $P$}\},
\end{align}
so the proposition above says that 
$C^\oo_{\adm}(P,V)\subset C^\oo_A(P,V)$.
It follows that for $h\in C^\oo_A(P,V)$ there exists a unique $\ga_h\in\gl(V)$ such that 
\begin{equation}\label{dis}
L_{\dot\rh(v)}h=-\ga_h(v),\quad\forall v\in V.
 \end{equation}


\paragraph{$\om$-Equivariant diffeomorphisms.}

There is a natural $(\Diff(M)\x\GL(V))$--action on the vector space $\Om^2(M,V)$ of $V$--valued 2--forms on $M$:
$$
(\ph,u)\cdot\om=u\cdot((\ph^{-1})^*\om)
$$
with infinitesimal action of the Lie algebra $\X(M)\x\gl(V)$
with Lie bracket (\ref{lb}) given by
$$
(\et,\ga)\cdot\om=-L_\et\om+\ga\cdot\om.
$$
The isotropy group of a closed 2--form $\om\in\Om^2(M,V)$ is the group of {\it $\om$--equivariant diffeomorphisms} 
\begin{equation*}
\Diff^{\eq}(M,\om):=\{(\ph, u)\in\Diff(M)\x\GL(V):\ph^*\om= u\cdot\om\}.
\end{equation*}
The isotropy Lie algebra coincides with the Lie algebra of 
$\om$--equivariant vector fields $\X^{\eq}(M,\om)$.

The group of $\om$--equivariant diffeomorphisms contains the group of $\om$--invariant diffeomorphisms as a subgroup
via the injective homomorphism 
\begin{equation*}
i:\Diff(M,\om)\to\Diff^{\eq}(M,\om),\quad i(\ph)=(\ph,1_V).
\end{equation*}
The restriction 
of the second component $ u\in\GL(V)$ of the $\om$--equivariant diffeomorphism $(\ph,u)$ to $V_0\subseteq V$ (defined in Remark \ref{vzero}) is determined by its first component $\ph\in\Diff(M)$. 
When $V_0=V$, then one identifies $\Diff^{\eq}(M,\om)$ with its projection on the first factor, the 
group $\{\ph\in\Diff(M):\exists u\in\GL(V)\text{ s.t. }\ph^*\om= u\cdot\om\}$.

\medskip
A curve $\ph$ in $\Diff(M)$ 
is called a {\it smooth curve} if the
corresponding map $(t,x)\mapsto (\ph(t)(x),\ph(t)^{-1}(x))$
in $M\x M$ is smooth.
Similarly a curve $u$ in $\GL(V)$ is smooth if the map $(t,v)\mapsto (u(t)(v),u(t)^{-1}(v))$ in $V\x V$ is smooth.
Let $\Diff^{\eq}(M,\om)_0$ be the normal subgroup of those elements in
$\Diff^{\eq}(M,\om)$ which can be connected to the identity by a smooth
curve in $\Diff^{\eq}(M,\om)\subset\Diff(M)\x\GL(V)$.

\begin{remark}\label{triv}
The second projection $(\ph, u)\in\Diff^{\eq}(M,\om)\mapsto u\in\GL(V)$ 
is a group homomorphism, 
so $V$ becomes a natural $\Diff^{\eq}(M,\om)$--module. 
The abelian group $A=V/\Ga$ is a $\Diff^{\eq}(M,\om)_0$--module too.
Indeed, for any 2--cycle $\si$ in $M$ and for any $(\ph, u)\in\Diff^{\eq}(M,\om)_0$,
the 2--cycles $\si$ and $\ph(\si)$ are homologous, so
$u(\int_\si\om)=\int_\si\ph^*\om=\int_\si\om$ and $u$
fixes the elements of the period group $\Ga$. In particular $u$ descends to a group automorphism 
$\bar u$ of $A$, and
the $\Diff^{\eq}(M)$--action on $V$ descends to an action on the abelian group $A$. If $V$ is generated by $\Ga$, then $u=1_V$ and $\ph$ is $\om$--invariant for all $\om$--equivariant diffeomorphisms $(\ph,u)$.
\end{remark}

\begin{proposition}\label{banach}
The following equivalences
hold for smooth paths $\ph_t$ in $\Diff(M)$ and $u_t$ in $\GL(V)$ starting at the identity:
$(\ph_t, u_t)\in\Diff^{\eq}(M,\om)$$\Leftrightarrow$
$(\de^l\ph_t,\de^l u_t)\in\X^{\eq}(M,\om)$
$\Leftrightarrow$$(\de^r\ph_t,\de^r u_t)$
$\in\X^{\eq}(M,\om)$.
\end{proposition}

This follows from Remark \ref{rele} in the Appendix.
In particular if the flow of an $\om$--equivariant vector field $(\et,\ga)$ exists, then it
consists of $\om$--equivariant diffeomorphisms.

\begin{lemma}\label{loop}
For any loop $\ell$ in $P$ and any smooth path of $\om$--equivariant 
diffeomorphisms $(\ph_t, u_t)$ starting at the identity,
we define the 2--chain $\si$ swept out by the loop $q\o\ell$ in $M$ under the isotopy $\ph_t$,
\ie $\si(t,s)=\ph_t(q(\ell(s)))$, $t,s\in [0,1]$. Then
\begin{equation*}
\int_\ell\int_0^1 u_t\cdot(q^*i_{\de^l\ph_t}\om
-\de^l u_t\cdot\th)dt
=\int_\si\om- u\cdot\int_\ell\th+\int_\ell\th.
\end{equation*}
\end{lemma}

\begin{proof}
Using $\ph_t^*\om=u_t\cdot\om$, we compute
\begin{align*}
\int_0^1\Big( u_t\cdot&\int_\ell q^*i_{\de^l\ph_t}\om\Big) dt=
\int_0^1 u_t\cdot \Big(\int_0^1(q^*i_{\de^l\ph_t}\om)(\dot \ell(s))ds\Big)dt\\
&=\int_0^1\int_0^1 u_{t}\cdot\om(\de^l\ph_t(q(\ell(s))),
Tq.\dot \ell(s))dsdt\\
&=\int_0^1\int_0^1\om(\dot\ph_t(q(\ell(s))),T\ph_t.Tq.\dot \ell(s))dsdt
=\int_\si\om
\end{align*}
and the result follows.
\end{proof}

\paragraph{Flux cocycle.}

The quotient space $H^1(P, V)/H^1(P,\Ga)$ receives a natural $\Diff^{\eq}(M,\om)_0$--module structure
from the action on the range $V$. Indeed, as we have seen
in Remark \ref{triv},
$\Diff^{\eq}(M,\om)_0$ acts trivially on $\Ga\subset V$,
hence it acts trivially on $H^1(P,\Ga)\subset H^1(P,V)$.

The map
\begin{align}\label{fluxeq}
&\Flux^{\eq}:{\Diff}^{\eq}(M,\om)_0\to H^1(P, V)/H^1(P,\Ga)\nonumber\\
\Flux^{\eq}(\ph,u)&=\int_0^1 u_t\cdot\flux^{\eq}(\de^l\ph_t,\de^lu_t)dt\mod H^1(P,\Ga)\nonumber\\
&=\Big[\int_0^1 u_t\cdot(q^*i_{\de^l\ph_t}\om
-\de^l u_t\cdot\th) dt\Big]\mod H^1(P,\Ga),
\end{align}
for any piecewise smooth path of $\om$--equivariant 
diffeomorphisms $(\ph_t, u_t)$ from the identity to $(\ph, u)$,
is a well defined group 1--cocycle,
called {\it the flux 1--cocycle}
associated to the closed vector valued form $\om$
with discrete period group $\Ga$.

The map $\Flux^{\eq}$ is well defined because Lemma \ref{loop} implies that for a loop $(\ph_t,u_t)$ of $\om$--equivariant diffeomorphisms based at the identity, the integral over a loop $\ell$ in $P$
of the 1--form $\int_0^1 u_t\cdot(q^*i_{\de^l\ph_t}\om
-\de^l u_t\cdot\th)dt$ 
is the integral of $\om$ over a 2--cycle $\si$, hence 
it belongs to the group $\Ga$ of periods of $\om$. 
The 1--cocycle condition for $\Flux^{\eq}$ is verified as in Proposition \ref{1cocy}
from the Appendix.

\begin{remark}
The group $\Diff_{\ex}^{\eq}(M,\om)=(\Ker \Flux^{\eq})_0$ is called the group 
of {\it equi-hamil\-tonian diffeomorphisms}.
The flux 1--cocycle $\Flux^{\eq}$ and the flux homomorphism $\Flux$ are related by
$\Flux^{\eq}\o i=q^*\o\Flux$, 
hence $i$ descends to an injective homomorphism $\Diff_{\ex}(M,\om))\to\Diff_{\ex}^{\eq}(M,\om)$.
\end{remark}

The next proposition follows from Remark \ref{refinement}
in the Appendix.

\begin{proposition}\label{eqex}
For any piecewise smooth path  $(\ph_t, u_t)$ of $\om$--equivariant 
diffeomorphisms, we have
$(\ph_t, u_t)\in\Diff_{\ex}^{\eq}(M,\om)$$\Leftrightarrow$
$(\de^l\ph_t,\de^l u_t)\in\X_{\ex}^{\eq}(M,\om)$$\Leftrightarrow$
$(\de^r\ph_t,\de^r u_t)\in\X_{\ex}^{\eq}(M,\om)$.
\end{proposition}


\begin{remark}\label{deal}
In the special case when $\om=d\al$ for an
$\al\in\Om^1(M, V)$ (in particular the period group $\Ga$ is trivial), the flux homomorphism is given by $\Flux(\ph)=[\ph^*\al-\al]\in H^1(M,V)$.
To compute the flux 1--cocycle $\Flux^{\eq}$ in this case,
let $P=M\x V\stackrel{q}{\to}M$ be the trivial $V$-bundle
with principal connection 1--form $\th=q^*\al+\th_V$ and curvature $\om$, 
where $\th_V=\de^l(1_V)\in\Om^1(V,V)$ stands for the Maurer-Cartan form on $ V$. We get
\begin{gather*}
\Flux^{\eq}:\Diff^{\eq}(M,\om)\to H^1(P, V),\quad
\Flux^{\eq}(\ph, u)=q^*[\ph^*\al- u\cdot\al].
\end{gather*}
Indeed, let $(\ph_t, u_t)$ be 
a path of $\om$--equivariant diffeomorphisms joining the identity and $(\ph,u)$. From $\tfrac{d}{dt}[\ph_t^*\al]=[u_t\cdot i_{\de^l\ph_t}\om]$ we obtain $\tfrac{d}{dt}[q^*\ph_t^*\al- u_t\cdot\th]=
u_t\cdot[q^*i_{\de^l\ph_t}\om-\de^l u_t\cdot\th]$.
Integrating this cohomology class from 0 to 1 gives 
the expression of the flux cocycle
$\Flux^{\eq}(\ph, u)=q^*[\ph^*\al- u\cdot\al]-[u\cdot\th_V-\th_V]=q^*[\ph^*\al- u\cdot\al]$,
because $\th_V$ is an exact 1--form on $P$.
\end{remark}


\section{Prequantization abelian extension}\label{s4}

Appropriate prequantization procedures have been developed for symplectic, presymplectic, Poisson

and Dirac manifolds \cite{WZ05}. In this section we suggest
a prequantization procedure for 
a closed vector valued 2--form in the equivariant setting.

\paragraph{Infinitesimal equi-quantomorphisms.}
The space $C^\oo_A(P,V)$ of almost $A$--invariant functions
defined in (\ref{capv})
endowed with the Lie bracket 
\begin{equation}\label{bra}
[h_1,h_2]=\ga_{h_2}\o h_1-\ga_{h_1}\o h_2,
\end{equation}
is a Lie subalgebra of $C^\oo(P,V)$ with Lie bracket (\ref{[]}),
since from \eqref{dis} follows that $L_{\dot\rh(h_1)}h_2=-\ga_{h_2}\o h_1$ for 
$h_1,h_2\in C^\oo_A(P,V)$.

\begin{align}\label{hieq}
\X^{\eq}(P)^A&=\{\xi\in\X(P):\exists\ga\in\gl(V)\text{ s.t. }\forall v\in V, L_{\dot\rh(v)}\xi=\dot\rh(\ga(v))\}
\end{align}
is the Lie algebra of {\it almost $A$--invariant vector fields}.
It can be characterized as
\begin{equation}\label{xepa}
\X^{\eq}(P)^A=
\{\xi\in\X(P):\exists\bar\ga:A\to V\text{ s.t. }\forall a\in A,
\rh(a)^*\xi-\xi=\dot\rh(\bar\ga(a))\},
\end{equation}
because of the identity
$$
\tfrac{d}{dt}(\rh(\exp tv)^*\xi-\xi-\dot\rh(\bar\ga(\exp tv)))
=\rh(\exp tv)^*(L_{\dot\rh(v)}\xi-\dot\rh(\ga(v))),
$$
where $\bar\ga$ is a group homomorphism with $\bar\ga\o\exp=\ga$. 
In particular every almost $A$--invariant vector field is projectable.

\begin{proposition}
When restricting the extension \eqref{proj} to the Lie algebra $\X^{\eq}(P)^A$
of almost $A$--invariant vector fields, one obtains
a new Lie algebra extension
\begin{equation}\label{infadm}
0\to C_A^\oo(P, V)\stackrel{\dot\rh}{\to}\X^{\eq}(P)^A \stackrel{q_*}{\to}\X(M)\to 0,
\end{equation}
with $C^\oo_A(P,V)$ the space of almost $A$--invariant functions defined in \eqref{capv}.
\end{proposition}

\begin{proof}
For an arbitrary $h\in C^\oo(P,V)$, the necessary and sufficient condition for
the vertical vector field $\dot\rh(h)$ to be almost $A$--invariant
is $L_{\dot\rh(v)}\dot\rh(h)=\dot\rh(\ga(v))$. But we know from (\ref{[]}) that $L_{\dot\rh(v)}\dot\rh(h)=-\dot\rh(L_{\dot\rh(v)}h)$, so the condition above becomes $L_{\dot\rh(v)}h=-\ga(v)$,
which means $h\in C^\oo_A(P,V)$.
\end{proof}

To pass to an abelian extension by $V$, we have 
to consider the Lie algebra 
\begin{equation*}
\X^{\eq}(P,\th)^A=\{\xi\in\X_{\proj}(P):\exists\ga_\xi\in\gl( V)
\text{ s.t. }L_\xi\th=\ga_\xi\cdot\th\}
\end{equation*}
of {\it infinitesimal equi-quantomorphisms}.
The linear map $\ga_\xi\in\gl(V)$ is determined by $\xi$ because
$\ga_\xi(v)=\ga_\xi(i_{\dot\rh(v)}\th)=i_{\dot\rh(v)}L_\xi\th$.
The Lie algebra $\X^{\eq}(P,\th)^A$ contains as a
Lie subalgebra the Lie algebra $\X(P,\th)^A$ of 
infinitesimal quantomorphisms.

An equivariant version of Remark \ref{corr} holds.

\begin{proposition}\label{pro}
The linear map
\begin{equation}\label{cor}
\xi\in\X^{\eq}(P,\th)^A\mapsto -i_\xi\th\in C_{\adm}^\oo(P,V).
\end{equation}
is surjective with kernel $\Ga(\Ker\om)^{\hor}$.
\end{proposition}

\begin{proof}
Given $h\in C^\oo_{\adm}(P,V)$, let $(\et,\ga)$ be
an equi-hamiltonian vector field with equi-hamiltonian function $h$
and let $\xi:=\et^{\hor}-\dot\rh(h)$. 
Then $\xi\in\X^{\eq}(P,\th)^A$ because
$L_{\xi}\th=L_{\et^{\hor}}\th-L_{\dot\rh(h)}\th=q^*i_{\et}\om-dh
=\ga\cdot\th$.
The linear correspondence (\ref{cor}) is surjective since $-i_{\xi}\th= h$.

Let $\xi\in\X^{\eq}(P,\th)^A$ be an element in the kernel of (\ref{cor}).
Then $\xi$ is a horizontal lift: there exists $\et\in\X(M)$
such that $\xi=\et^{\hor}$. But $\xi$ is an infinitesimal quantomorphism,
so $\ga_\xi\cdot\th=L_{\et^{\hor}}\th=i_{\et^{\hor}}d\th=q^*i_\et\om$.
Evaluating on vertical vectors we get $\ga_\xi=0$, so $i_{\et}\om=0$, which means $\et\in\Ga(\Ker\om)$.
\end{proof}

The next proposition shows the inclusion $\X^{\eq}(P,\th)^A\subset\X^{\eq}(P)^A $.

\begin{proposition}\label{incl}
Any infinitesimal equi-quantomorphism with $L_\xi\th=\ga_\xi\cdot\th$ 
satisfies
$L_{\dot\rh(v)}\xi=\dot\rh(\ga_\xi(v))$ 
for all $v\in V$.
\end{proposition}

\begin{proof}
The infinitesimal equi-quantomorphism $\xi$ is projectable 
and $\dot\rh(v)$ is vertical, so $L_{\dot\rh(v)}\xi=[\dot\rh(v),\xi]$ is also vertical.
A short computation using $L_\xi\th=\ga_\xi\cdot\th$
gives $i_{[\dot\rh(v),\xi]}\th=\ga_\xi(v)$.
These two facts imply $L_{\dot\rh(v)}\xi=\dot\rh(\ga_\xi(v))$.
\end{proof}

The infinitesimal equi-quantomorphism $\xi$ determines an
$\om$--equivariant vector field $(q_*\xi,\ga_\xi)$ on $M$.
Moreover $(q_*\xi,\ga_\xi)$ is an equi-hamiltonian vector field
for the equi-hamiltonian function 
$h=-i_\xi\th\in C^\oo_{\adm}(P,V)$, because
$$
dh=-di_\xi\th=i_\xi d\th-L_\xi\th=q^*i_{q_*\xi}\om-\ga_\xi\cdot\th.
$$

\begin{theorem}\label{infini}
The Lie algebra of infinitesimal equi-quantomorphisms is an abel\-ian
extension of the Lie algebra $\X_{\ex}^{\eq}(M,\om)$ of equi-hamiltonian vector fields by 
the natural $\X_{\ex}^{\eq}(M,\om)$--module $ V$.
An abelian Lie algebra 2--cocycle on $\X_{\ex}^{\eq}(M,\om)$ 
defining this abelian extension
is $((\et_1,\ga_1),(\et_2,\ga_2))\mapsto -\om(\et_1,\et_2)(x_0)$, for any fixed element $x_0\in M$.
\end{theorem}

\begin{proof}
We show that the following sequence of Lie algebras is exact:
\begin{equation}\label{algebraextension}
0\to V\stackrel{\dot\rh}{\to}\X^{\eq}(P,\th)^A
\stackrel{p}{\to}\X^{\eq}_{\ex}(M,\om)\to 0,
\end{equation}
where $p(\xi)=(q_*\xi,\ga_\xi)$.
The injectivity of $\dot\rh$ is clear.
For the surjectivity of $p$
we consider an equi-hamiltonian vector field $(\et,\ga)$ with equi-hamiltonian
function $h\in C^\oo(P, V)$.
Then the infinitesimal equi-quantomorphism
$\xi=\et^{\hor}-\dot\rh(h)$ (from the proof of Proposition \ref{pro}) projects to $(\et,\ga)$.

The inclusion $\dot\rh(V)\subseteq\Ker p$ follows from $p\o\dot\rh=0$.
For the reversed inclusion let 
$\xi\in\Ker p\subset\X^{\eq}(P,\th)^A$. 
Then $L_\xi\th=0$ and $\xi=\dot\rh(h)$ for some $h\in C^\oo(P, V)$.
From $L_{\dot\rh(h)}\th=dh$ and from the connectedness of $P$  follows that $\xi\in\dot\rh( V)$.
The induced action of $\X^{\eq}_{\ex}(M,\om)$ on $ V$ is 
the natural one because from Proposition \ref{incl}
we get $[\dot\rh(v),\xi]=\dot\rh(\ga_\xi(v))$ for all $v\in V$ and $\xi\in\X^{\eq}(P,\th)^A$.

We determine the 2--cocycle defined with the linear section $s$
of (\ref{algebraextension}) given by 
$s(\et,\ga)=\et^{\hor}-\dot\rh(h)$, where $h$ is the unique equi-hamiltonian function of the equi-hamiltonian vector field $(\et,\ga)$
vanishing at a fixed point $y_0\in q^{-1}(x_0)$.
First we observe that given
the equi-hamiltonian vector fields
$(\et_1,\ga_1)$ and $(\et_2,\ga_2)$ with equi-hamiltonian 
functions $h_1$ and $h_2$ vanishing at $y_0$,
the equi-hamiltonian function vanishing at $y_0$ for
the bracket $(-[\et_1,\et_2],\ga_1\ga_2-\ga_2\ga_1)$ is
$\ga_1\o h_2-\ga_2\o h_1+q^*\om(\et_1,\et_2)-\om(\et_1,\et_2)(x_0)$.
Indeed,
\begin{align*}
d(\ga_1\o h_2-\ga_2\o h_1&+q^*\om(\et_1,\et_2))\\
&\stackrel{\eqref{defham}}{=}\ga_1\cdot(q^*i_{\et_2}\om-\ga_2\cdot\th)
-\ga_2\cdot(q^*i_{\et_1}\om-\ga_1\cdot\th)
+q^*di_{\et_2}i_{\et_1}\om\\
&=q^*(i_{\et_2}L_{\et_1}\om-i_{\et_1}L_{\et_2}\om+di_{\et_2}i_{\et_1}\om)-(\ga_1\o\ga_2-\ga_2\o\ga_1)\cdot\th\\
&=-q^*i_{[\et_1,\et_2]}\om-(\ga_1\o\ga_2-\ga_2\o\ga_1)\cdot\th.
\end{align*}
Since $\om$ is the curvature of the principal bundle $P$, the identity
$[\et_1^{\hor},\et_2^{\hor}]-[\et_1,\et_2]^{\hor}=
-\dot\rh(q^*\om(\et_1,\et_2))$ holds.
The equi-hamiltonian functions $h_1$ and $h_2$ belong to $C^\oo_A(P,V)$, so by
(\ref{bra}) 
we have
$[\dot\rh(h_1),\dot\rh(h_2)]=\dot\rh(\ga_1\o h_2-\ga_2\o h_1)$.
Using also the fact that
$$[\dot\rh(h_1),\et_2^{hor}]=
-\dot\rh(L_{\et_2^{hor}}h_1)\stackrel{\eqref{defham}}{=}-\dot\rh(q^*\om(\et_1,\et_2)), 
$$
we compute
\begin{align*}
[s(\et_1,\ga_1),s(\et_2,\ga_2)]&-s([(\et_1,\ga_1),(\et_2,\ga_2)])
=-[\et_1^{\hor}-\dot\rh(h_1),\et_2^{\hor}-\dot\rh(h_2)]\\
&+[\et_1,\et_2]^{\hor}
-\dot\rh(\ga_2\o h_1-\ga_1\o h_2-q^*\om(\et_1,\et_2)
+\om(\et_1,\et_2)(x_0))\\
&=[\dot\rh(h_1),\et_2^{\hor}]-[\dot\rh(h_2),\et_1^{\hor}]
-[\dot\rh(h_1),\dot\rh(h_2)]\\
&-\dot\rh(\ga_2\o h_1)
+\dot\rh(\ga_1\o h_2)+2\dot\rh(q^*\om(\et_1,\et_2))
-\dot\rh(\om(\et_1,\et_2)(x_0))\\
&=-\dot\rh(\om(\et_1,\et_2)(x_0)),
\end{align*}
thus obtaining a Lie algebra 2--cocycle for the abelian extension (\ref{algebraextension}).
\end{proof} 

\begin{remark}\label{like}
Under the assumption that the closure $V_0$ of
the image of $\om:TM\x_M TM\to V$ is the whole $V$, 
the Lie algebra of equi-hamiltonian vector fields can  be identified with its projection on the first factor:
$$\X_{\ex}^{\eq}(M,\om)=\{\et\in\X(M):\exists \ga\in\gl(V)\text{ s.t. }q^*i_\et\om-\ga\cdot\th\text{ exact}\}. 
$$
In this case the abelian extension
(\ref{algebraextension}) can be seen as a restriction of (\ref{infadm}).
\end{remark}


\paragraph{Prequantization abelian extension.}
We saw in Sect.~\ref{s2} that the gauge extension \eqref{a} (of $A$-equivariant diffeomorphisms of $P$)
is an abelian group extension which contains the prequantization central extension.
In this paragraph we construct the prequantization abelian extension
together with a non-abelian extension (of almost $A$-equivariant diffeomorphisms of $P$)
containing it.

The group of {\it almost $A$--equivariant diffeomorphisms} of $P$ is 
$$
\Diff_A(P)=\{\ps\in\Diff(P):\exists \bar u_\ps\in\Aut(A)\text{ s.t. }
\ps\o\rh(a)=\rh(\bar u_\ps(a))\o\ps,\forall a\in A\},
$$
where $\Aut(A)$ denotes the group of automorphisms of the abelian group $A$.
One can describe $\Diff_A(P)$ as the group of those diffeomorphisms $\ps$ of $P$ such that the vertical vector fields $\dot\rh(v)$ and $\dot\rh(u_\ps(v))$ are $\ps$-related.
It is a subgroup of the group $\Diff_{\proj}(P)$ of projectable diffeomorphisms
and it contains the group $\Diff(P)^A$ of equivariant diffeomorphisms of $P$ as a subgroup.

The fiber preserving diffeomorphism $\rh(f)$ for $f\in C^\oo(P,A)$ is almost $A$--equi\-variant
if and only if $\rh(f)\o\rh(a)=\rh(\bar u_{\rh(f)}(a))\o\rh(f)$,
which can be written as $f(y)^{-1}f(\rh(y,a))=a^{-1}\bar u_{\rh(f)}(a)$ for all $y\in P$.
We define the set of {\it almost $A$--invariant maps}
\begin{align*}
C^\oo_{A}(P,A)
=\{f\in C^\oo(P,A): \forall a\in A, f^{-1}(f\o\rh(a)) \text{ constant on $P$}\}.
\end{align*}
For an almost $A$--invariant function $h\in C^\oo_A(P,V)$,
$\exp\o h$ is an almost $A$--invariant map with $\bar u_f=\exp\o\bar\ga_h$.
The existence of a unique $\bar u_f\in\Aut(A)$ such that  
$f(y)^{-1}f(\rh(y,a))=a^{-1}\bar u_f(a)$ for all $a\in A$ and $y\in P$
follows easily.
We observe that $\rh(f)\in\Diff_A(P)$ if and only if 
$f$ is an almost $A$--invariant map with $\bar u_f=\bar u_{\rh(f)}$. 

With the group multiplication
$$
(f_1\cdot f_2)(y)=f_1(y)\bar u_{f_1}(f_2(y)),
$$
$C^\oo_{A}(P,A)$  becomes a
subgroup of $C^\oo(P,A)_{\inv}$ with multiplication (\ref{dot}). 
The abelian group $C^\oo(M,A)$, identified with the group of $A$--invariant maps $f:P\to A$, is a subgroup of $C_A^\oo(P,A)$ (in this case $\bar u_f=1_A$).

Restricting (\ref{star}) to almost $A$-equivariant diffeomorphisms, 
we obtain a non-abelian group extension with infinitesimal version (\ref{infadm}): 
\begin{equation}\label{starstar}
1\to C_{A}^\oo(P,A)\stackrel{\rh}{\to}\Diff_A(P)
\stackrel{q_*}{\to}\Diff(M)_{[P]}\to 1.
\end{equation}
This is an enlarged version of the gauge extension \eqref{a}.


We define the group of {\it equi-quanto\-morphisms} as
\begin{equation*}
\Diff^{\eq}(P,\th)^A=\{\ps\in\Diff_{\proj}(P):\exists 
 u_\ps\in\GL( V)\text{ s.t. }\ps^*\th= u_\ps\cdot\th\}.
\end{equation*}
It contains the group $\Diff(P,\th)^A$ of quantomorphisms as a subgroup.
The linear isomorphism $ u_\ps$ is uniquely determined by $\ps\in\Diff^{\eq}(P,\th)^A$.
Let $\ph$ denote the diffeomorphisms of $M$ with $q\o\ps=\ph\o q$. By differentiating
the relation $\ps^*\th= u_\ps\cdot\th$ we get $\ph^*\om= u_\ps\cdot\om$,
hence $(\ph, u_\ps)$ is an $\om$--equivariant diffeomorphism.

The proof of the next proposition follows from Remark \ref{rele} in the Appendix.

\begin{proposition}\label{equi}
For a smooth curve $\ps_t$ in $\Diff(P)$ starting at the identity we have
$\ps_t\in\Diff^{\eq}(P,\th)^A\Leftrightarrow\de^l\ps_t\in\X^{\eq}(P,\th)^A
\Leftrightarrow\de^r\ps_t\in\X^{\eq}(P,\th)^A$. 
\end{proposition}

In particular if the flow of an infinitesimal equi-quantomorphism exists, then it consists of
equi-quantomorphisms.

\begin{proposition}\label{subs}
The group $\Diff_A(P)$ of almost $A$--equivariant   diffeomorphisms contains the group $\Diff^{\eq}(P,\th)^A$ of equi-quantomorphisms as a subgroup. 
More precisely, the deviation from $A$--equivariance of an
equi-quantomorphism $\ps$ is measured by the isomorphism $u_\ps\in\GL(V)$. 
\end{proposition}

\begin{proof}
Let $\ps$ be an equi-quantomorphism with $\ps^*\th=u_\ps\cdot\th$. Then $\ps^{-1}$ is an equi-quantomorphism too, with $(\ps^{-1})^*\th=(u_\ps)^{-1}\cdot\th$. For any $a\in A$, the diffeomorphism 
$\ps\o\rh(a)\o\ps^{-1}$ is fiber preserving and $\th$--invariant,
in particular it is of the form $\rh(f)$ with $f\in C^\oo(P,A)_{\inv}$ depending on $a$.
From Proposition \ref{appendix} in the Appendix, $0=\rh(f)^*\th-\th=\de^l(f)$, so that $f$ is a constant denoted $\bar u(a)\in A$.
We obtain that $\ps\o\rh(a)=\rh(\bar u(a))\o\ps$,
so $\bar u\in\Aut(A)$.
The infinitesimal version of this identity is $T\ps\o\dot\rh(v)=\dot\rh(u(v))\o\ps$,
where $u\in\GL(V)$ with $\bar u\o\exp=\exp\o u$.

It remains to be shown that $u=u_\ps$.
This follows from the above mentioned fact that $\dot\rh(v)$ and $\dot\rh(u(v))$ are $\ps$-related:
$u_\ps(v)=i_{\dot\rh(v)}(\ps^*\th)=\ps^*(i_{\dot\rh(u(v))}\th)=u(v)$ for all $v\in V$. 
\end{proof}

The group $\Hol^{\eq}(M,\om)$ of diffeomorphisms preserving the holonomy up to a group automorphism of $A$,
called the group of {\it almost holonomy preserving diffeomorphisms} is
\begin{equation*}
\Hol^{\eq}(M,\om)=\{(\ph,\bar u)\in\Diff(M)\x\Aut(A):\forall\ell\in C^\oo(S^1,M),h(\ph\o\ell)=\bar u(h(\ell))\}.
\end{equation*}
The group $\Hol^{\eq}(M,\om)$ acts in a natural way on the abelian group $A$. 
Adapting the idea of the proof of Theorem 2.7 in \cite{NV03} to the equivariant setting, we will show 
in the appendix that
the group of equi-quantomorphisms is an abelian extension of the group of almost holonomy preserving diffeomorphisms.

\begin{lemma}
Given $(\ph,\bar u)\in\Hol^{\eq}(M,\om)$ and, for a fixed $x_0\in M$,
a bijection $\ps_{x_0}:q^{-1}(x_0)\to
q^{-1}(\ph(x_0))$ satisfying 
\begin{equation}\label{defin}
\ps_{x_0}\o\rh(a)=\rh(\bar u(a))\o\ps_{x_0},\quad\forall a\in A,
\end{equation}
there exists a unique equi-quantomorphism $\ps$ of $P$ extending $\ps_{x_0}$ and descending to the diffeomorphism $\ph$ of $M$. 
\end{lemma}

\begin{proof}
Let $\Pt(c):q^{-1}(x_0)\to q^{-1}(x)$ denote the parallel transport map along a curve $c$ from $x_0$ to $x$ in $M$.
It defines a map 
\begin{equation}\label{stea}
\ps_x=\Pt(\ph\o c)\o\ps_{x_0}\o\Pt(c)^{-1}:q^{-1}(x)\to q^{-1}(\ph(x))
\end{equation}
which does not depend on the choice of $c$, because for every loop $\ell$ at $x_0$
$$
\Pt(\ph\o\ell)\o\ps_{x_0}\o\Pt(\ell)^{-1}
=\rh(\bar u(h(\ell)))\o\ps_{x_0}\o\rh(h(\ell))^{-1}
\stackrel{\eqref{defin}}{=}\ps_{x_0}.
$$
The maps $\ps_x$, $x\in M$, glue to a diffeomorphism $\ps$ of $P$ which satisfies $\ps\o\rh(a)=\rh(\bar u(a))\o\ps$ for all $a\in A$.
Its infinitesimal version is: $\dot\rh(v)$ and $\dot\rh(u(v))$ are $\ps$-related, 
where $u\in\GL(V)$ is given by $\exp\o u=\bar u\o\exp$.

The tangent map $T\ps:TP\to TP$ maps horizontal vectors to horizontal vectors because for any horizontal lift $c^{\hor}$ of the curve $c$, the curve $\ps\o c^{\hor}$ is the horizontal lift starting at $\ps(c^{\hor}(0))$ of the curve $\ph\o c$.
Indeed, 
$\Pt(c)^{-1}(c^{\hor}(t))=c^{\hor}(0)$, so by \eqref{stea} we obtain 
$(\ps\o c^{\hor})(t)=\Pt(\ph\o c|_{[0,t]})(\ps(c^{\hor}(0))$.
Now one can show that $\ps^*\th=u\cdot\th$:
$$
(\ps^*\th)(\et^{\hor}+\dot\rh(v))=\th(T\ps.\et^{\hor})+\th(\dot\rh(u(v)))=u(v)=(u\cdot\th)(\et^{\hor}+\dot\rh(v))
$$
for all $\et\in\X(M)$ and $v\in V$.
\end{proof}

\begin{theorem}\label{kostant}
The group $\Diff^{\eq}(P,\th)^A$
of equi-quantomorphisms is an abelian extension of the group $\Hol^{\eq}(M,\om)$ of almost holonomy preserving diffeomorphisms by the natural 
$\Hol^{\eq}(M,\om)$--module $A$, \ie 
\begin{equation}\label{alho}
1\to A\stackrel{\rh}{\to}\Diff^{\eq}(P,\th)^A
\stackrel{p}{\to}\Hol^{\eq}(M,\om)\to 1
\end{equation} 
is an exact sequence of groups.
\end{theorem}

\begin{proof}
For $\ps\in\Diff^{\eq}(P,\th)^A$ with $\ph\o q=q\o\ps$ and $\ps^*\th=u_\ps\cdot\th$, we define $p(\ps)=(q_*\ps,\bar u_\ps)$,
where $\bar u_\ps\in\Aut(A)$ with $\bar u_\ps\o\exp=\exp\o u_\ps$.
In particular $\ps$ is almost $A$--equivariant: $\ps\o\rh(a)=\rh(\bar u_\ps(a))\o\ps$ by Proposition \ref{subs}.

We verify that $h(\ph\o\ell)=\bar u_\ps(h(\ell))$
for any $\ell\in C^\oo(S^1,M)$, 
showing that $\ph$ is an almost holonomy preserving diffeomorphism. 
First we observe that if $\ell^{\hor}$ is a horizontal lift of the loop $\ell$, then $\ps\o\ell^{\hor}$ is a horizontal lift of the loop $\ph\o\ell$:
$$
\th((\ps\o\ell^{\hor})'(t))=\th(T\ps.(\ell^{\hor})'(t))
=(\ps^*\th)((\ell^{\hor})'(t))=u_\ps(\th((\ell^{\hor})'(t))=0.
$$
Then the desired identity follows from the computation
\begin{align*}
\rh((\ps\o\ell^{\hor})(0),h(\ph\o\ell))&=(\ps\o\ell^{\hor})(1)
=\ps(\rh(\ell^{\hor}(0),h(\ell)))\\
&=\rh(\ps(\ell^{\hor}(0)),\bar u_\ps(h(\ell))).
\end{align*}
Thus $p:\Diff^{\eq}(P,\th)^A{\to}\Hol^{\eq}(M,\om)$ is well defined.
By the previous lemma it is also surjective, and \eqref{alho}
is an exact sequence of groups.
\end{proof}


Adapting the proof of Corollary 2.8 in \cite{NV03} to the equivariant setting, one shows that the identity component 
$\Diff^{\eq}(P,\th)^A_0$
of the group of equi-quanto\-morphisms is an abelian extension of the group $\Diff_{\ex}^{\eq}(M,\om)$ of 
equi-hamiltonian diffeomorphisms by the natural 
$\Diff_{\ex}^{\eq}(M,\om)$--module $A$:
\begin{equation}\label{four}
1\to A\stackrel{}{\to}\Diff^{\eq}(P,\th)^A_0
\stackrel{}{\to}\Diff_{\ex}^{\eq}(M,\om)\to 1.
\end{equation}
We call this the {\it prequantization abelian extension}.

\begin{remark}
Like in Remark \ref{like}, if $V=V_0$ the closure of
the image of $\om:TM\x_M TM\to V$, the group of equi-hamiltonian diffeomorphisms can be identified with its projection on the first factor. In this case both \eqref{alho} and \eqref{four} can be seen as a restriction of \eqref{starstar}.
\end{remark}


\section{Group extensions via prequantization}\label{s5}

The results from  Sect.~\ref{s2} concerning central Lie group extensions as pullback of the prequantization extension
are generalized in this section to obtain abelian Lie group extensions of $G$ associated to a $G$--equivariant 2--form $\om$.
Moreover, for simply connected $G$, all its abelian extensions can be obtained in this way.
 
We consider a smooth action $\la$ on $M$ and a linear action $b$ on $ V$
such that $\om$ is $G$--equivariant, \ie $\la_g^*\om=b_g\cdot\om$.
Let $\dot\la:\g\to\X(M)$ and $\dot b:\g\to\gl(V)$ denote the infinitesimal $\g$-actions.
The pair $(\la,b)$ is called an {\it equi-hamiltonian $G$--action} if the 1--form $q^*i_{\dot \la(X)}\om-\dot b(X)\cdot\th$ is exact for all $X\in\g$.

In this case $(\la_g,b_g)\in\Diff_{\ex}^{\eq}(M,\om)$ for each $g\in G$, by Proposition \ref{eqex}.
Then the pull-back of the prequantization abelian extension (\ref{four}) by $(\la,b)$: 
\begin{equation}\label{pullback}
\hat G=\{(g,\ps)\in G\x\Diff(P):q_*\ps=\la_g, \ps^*\th=b_g\cdot\th\},
\end{equation}
is 
an abelian group extension 
$1\to A\stackrel{i}{\to}\hat G\stackrel{p}{\to}G\to 1$ with $i(a)=(e,\rh(a))$ and $p(g,\ps)=g$. 
The induced $G$--module structure on $A$ comes from the linear action $b$ on $V$.

To show that $\hat G$ is a Lie group with a smooth action $\hat \la$  
by equi-quantomorphisms of $(P,\th)$, lifting the action $\la$, we use Lemma 3.2 and Lemma 3.3 from \cite{NV03} and we adapt the proof of Theorem 2.7 from \cite{NV03} to abelian Lie group extensions.
This is resumed in the next theorem.

\begin{theorem}\label{th2}
Given $\om\in\Om^2(M,V)$ a closed 2--form with discrete period group $\Ga$
and an equi-hamiltonian $G$--action $(\la,b)$, there is an abelian
Lie group extension $\hat G$ of $G$ by the $G$--module $A=V/\Ga$,
integrating the $V$--valued Lie algebra 2--cocycle on $\g$
\begin{equation}\label{dotdot}
(X,Y)\mapsto-\om(\dot \la(X),\dot \la(Y))(x_0),
\end{equation}
whose cohomology class does not
depend on the choice of the point $x_0\in M$.
There is also a $\hat G$--action $\hat\la$ on $P$, lifting the $G$--action $\la$, and 
such that $\th$ is $\hat G$--equivariant,
\ie $\hat \la_{\hat g}^*\th=b_{p(\hat g)}\cdot\th$ for any $\hat g\in\hat G$.
\end{theorem}

\begin{proof}
Let $\la^{x_0}:g\in G\mapsto \la_g(x_0)\in M$ be the orbit map,
and let $(\la^{x_0})^*P\to G$ be the pull-back of the $A$-bundle $P\to M$. 
From Lemma 3.2 in \cite{NV03} follows that an element $y_0\in P$
with $q(y_0)=x_0$ defines a bijection $(g,\ps)\in\hat G\mapsto (g,\ps(y_0))\in (\la^{x_0})^*P$, 
and the smooth manifold structure transported on $\hat G$ by this bijection 
does not depend on the choice of $x_0$.

We define the $\hat G$--action $\hat \la$ on $P$ by 
$$
\hat \la:((g,\ps),y)\in\hat G\x P\mapsto\ps(y)\in P.
$$ 
It lifts the $G$--action $\la:G\x M\to M$
and $\th$ is $\hat G$--equivariant because $\ps$ is an equi-quantomorphism. 
Its restriction to $A$ is 
$\hat \la_a=\rh(a)$
for $a\in A$.

In order to show that $\hat \la$ is smooth, 
we put product coordinates on $P$ and $\hat G$.
Any smooth local section $s_M:U\subset M\to P$ of $q$, with $s_M(x_0)=y_0$,
defines a local smooth section $s_G(g)=(g,\ps_g)\in\hat G$ of $p$,
for $g$ sufficiently close to the identity of $G$ in order to have
$\la_g(x_0)\in U$. Here $\ps_g$ is uniquely defined by the condition $\ps_g(y_0)=s_M(\la_g(x_0))$.

We fix an open neighbourhood $U_G$ of $e$ and a neighbourhood 
$U_M$ of $x$ diffeomorphic to an open convex subset of the modeling space,
such that $U_G\cdot U_M\subset U$. Then we define a 
function $f:U_G\x U_M\to A$ with $f(g,x_0)=1$ for all $g\in U_G$ by the relation 
\begin{equation}\label{auto}
\ps_g(s_M(x))=\rh(s_M(\la_g(x)),f(g,x)).
\end{equation}
In other words,
for $g\in U_G$ the expression of $\ps_g$ in product
coordinates is $\ps_g:(x,a)\in U_M\x A\mapsto(\la_g(x),f(g,x)\bar b_g(a))\in U\x A$, where $\bar b_g\in\Aut(A)$ with $\bar b_g\o\exp=\exp\o b_g$.
This can be deduced from
\begin{align*}
(\ps_g\o\rh(a))(s_M(x))
=(\rh({\bar b_g(a)})\o\ps_g)(s_M(x))
=\rh(f(g,x)\bar b_g(a))(s_M(\la_g(x))).
\end{align*}

The connection 1--form in product coordinates 
$U\x A$ is 
$\th=q^*\al+q_A^*\th^l_A$, with $q_A:q^{-1}(U)\to A$ the second projection, $\al=s_M^*\th\in\Om^1(U,V)$ and $\th^l_A=\de^l(1_A)\in\Om^1(A,V)$
the Maurer-Cartan form on $A$. In particular $d\al=\om$ on $U\subset M$. 
The condition $\ps_g^*\th=b_g\cdot\th$ implies
$\de^l(f_g)= \la_g^*\al -b_g\cdot\al $. Indeed,
\begin{align*}
0=\ps_g^*(q^*\al&+q_A^*\th^l_A)-b_g\cdot(q^*\al+q_A^*\th^l_A)
=q^*(\la_g^*\al-b_g\cdot\al)\\
&+\de^l(m\o(\bar b_g\o q_A,f_g\o q))-q_A^*\de^l(\bar b_g)
=q^*(\la_g^*\al-b_g\cdot\al-\de^l(f_g))
\end{align*}
because $q_A\o\ps_g=m\o(\bar b_g\o q_A,f_g\o q)$ for $m$ the group multiplication map of $A$ and
$\de^l(m\o(h_1,h_2))=\de^lh_1+\de^lh_2$ for $h_1,h_2:P\to A$.

Using the Poincar\'e Lemma applied to the closed 1--form 
$\la_g^*\al -b_g\cdot\al$ on the convex set $U_M$,
we obtain that $f$ is a smooth function.
Then
$$
(\rh(a)\o\ps_g)(\rh(s_M(x),a'))
=\rh(s_M(\la_g(x)),f(g,x)a\bar b_g(a'))
$$
assures the following expression of $\hat \la$ in product coordinates:
$$
\hat \la((g,a),(x,a'))=(\la_g(x),f(g,x)a\bar b_g(a')),
$$
thus showing the smoothness of $\hat \la$.

The $A$--valued local group 2--cocycle on $G$ corresponding to the section $s_G$
is 
\begin{equation}\label{name}
c(g_1,g_2)=f(g_1,\la_{g_2}(x_0)),
\end{equation}
because 
\begin{align*}
\ps_{g_1}\ps_{g_2}(y_0)&=\ps_{g_1}(s_M(\la_{g_2}(x_0)))
\stackrel{(\ref{auto})}{=}\rh(s_M(\la_{g_1g_2}(x_0)),f(g_1,\la_{g_2}(x_0)))\\
&=\rh(\ps_{g_1g_2}(y_0),f(g_1,\la_{g_2}(x_0))).
\end{align*}
This shows the smoothness of multiplication and inversion 
in an identity neighborhood in $\hat G$. That the left multiplications are smooth, follows from the fact that $\hat G$ acts by smooth maps on
the bundle $P$. Now Lemma 3.3 in \cite{NV03}
implies that $\hat G$ is a Lie group.

The corresponding Lie algebra extension $0\to V\to\hat\g\to\g\to 0$ is the pull-back by $(\dot \la,\dot b):\g\to\X_{\ex}^{\eq}(M,\om)$ of the abelian extension (\ref{algebraextension}),
hence a defining Lie algebra 2--cocycle is (\ref{dotdot}), the pull-back by $(\dot \la,\dot b)$ of the 2--cocycle from Theorem \ref{infini}.
\end{proof}

\begin{theorem}\label{th1}
Every abelian Lie group extension 
$1\to A=V/\Ga\to\hat G\stackrel{q}{\to} G\to 1$ of a simply connected regular Lie group $G$ can be obtained
as a pull-back of the prequantization abelian extension \eqref{four}. 
\end{theorem}

\begin{proof}
Let $0\to V\to\hat\g\stackrel{q}{\to}\g\to 0$ be the abelian Lie algebra extension corresponding to the given
Lie group extension, and let $s:\g\to\hat\g$ be a continuous section. We denote by $\si$ 
the $V$-valued 2--cocycle on $\g$ defined by $s$, and by $p:\hat\g\to V$ the associated projection, so
$\si(X,Y)=[s(X),s(Y)]-s([X,Y])$ and $p(\hat X)=\hat X-s(q(\hat X))$. 
We use the identification $\hat X=(X,v)$ for $X=q(\hat X)$ and $v=p(\hat X)$ 
(in this proof the same letter $q$ denotes both projections $\hat G\to G$ and $\hat\g\to\g$).

Let $\la$ be the left translation on $G$ and $b$ the linear $G$-action on $V$ 
induced by the $G$-module structure of $A$ determined by the abelian extension of $G$ by $A$,
namely $g\cdot a=\hat g a\hat g^{-1}$ for any $\hat g\in\hat G$ with $q(\hat g)=g$. 
With respect to these two $G$-actions, the $G$-equivariant 2--form 
$\om:=-\si^{\eq}\in\Om^2(G,V)$  is closed,
because $\si$ is a 2--cocycle. Since $\la_g^*\om=b_g\cdot\om$,
we obtain a group homomorphism $(\la,b):G\to\Diff^{\eq}(G,\om)$.
We wish to have a group homomorphism into $\Diff_{\ex}^{\eq}(G,\om)$,
and for this we need a principal connection on a principal bundle over $G$
with curvature $\om$.

The abelian extension $\hat G$ is a principal $A$-bundle over $G$.
Here we need the principal $A$-action on $\hat G$ to be given by left translations,
hence $\rh$ is a left action: $\rh(a)\hat g=a\hat g$. 
We show that the $\hat G$-equivariant 1-form $\th=p^{\eq}\in\Om^1(\hat G,V)$ (in the sense 
that $\hat\la_{\hat g}^*\th=b_g\cdot\th$ where $\hat\la$ denotes left translation on $\hat G$) 
with identity value $\th_e=p$,
is a principal connection. For this we only have to verify the identity 
$\th(\dot\rh(v))=v$ for all $v\in V$, because the $A$-invariance $\rh(a)^*\th=\th$ for all $a\in A$
is a special case of the $\hat G$-equivariance:
\begin{align*}
\th(\dot\rh(v))(\hat g)=\th\Big(\frac{d}{dt}\Big|_0(\exp tv)\hat g\Big)
=p^{\eq}(\hat g\cdot\Ad_{{\hat g}^{-1}}(v))
=b_g \cdot p(b_{g^{-1}}(v))=v.
\end{align*}
The curvature of the principal connection $\th$ is $\om$, 
because $p([(X_1,v_1),(X_2,v_2)])=\dot b(X_1)v_2-\dot b(X_2)v_1+\si(X_1,X_2)$
implies that the Chevalley-Eilenberg differential of the $V$-valued 1--cochain $p$ on $\hat\g$ 
is the pull-back of $-\si$ to $\hat\g$,
so $d\th=dp^{\eq}=-q^*\si^{\eq}=q^*\om$.

Next we show that the $G$-action $(\la,b)$ is equi-hamiltonian,
\ie the infinitesimal flux cocycle $\flux^{\eq}$ vanishes on 
$(\dot\la(Y),\dot b(Y))=(Y^r,\dot b(Y))\in\X^{\eq}(G,\om)$
for all $Y\in\g$, where $Y^r$ denotes the right invariant vector field on $G$ with identity value $Y^r(e)=Y$:
\begin{equation}\label{inter}
\flux^{\eq}(Y^r,\dot b(Y))=[q^*i_{Y^r}\om-\dot b(Y)\cdot\th]
=[q^*i_{Y^r}\si^{\eq}+\dot b(Y)\cdot p^{\eq}]\in H^1(\hat G,V).
\end{equation}
We will use  Proposition \ref{properties} from the Appendix to show that the 1--form 
$q^*i_{Y^r}\si^{\eq}+\dot b(Y)\cdot p^{\eq}$ on $\hat G$ is exact,
namely we will show that it is the equivariant 1-form coming from a Lie algebra 1--cocycle
which can be integrated to a group 1--cocycle. 

There is a group 1--cocycle $\ka:\hat G\to \Lin(\g,V)$ involved in the expression of the adjoint action 
on the abelian extension $\hat G$ \cite{N04}:
$$
\Ad_{\hat g}(Y,w)=(\Ad_gY,b_g(w)-\ka(\hat g)(\Ad_gY)).
$$
From the expression of the Lie bracket in $\hat\g$:
\[
\ad((X,v))(Y,w)=[(X,v),(Y,w)]=([X,Y],\dot b(X)w-\dot b(Y)v+\si(X,Y),
\]
we deduce that the Lie algebra 1--cocycle $\al:\hat\g\to \Lin(\g,V)$ corresponding to $\ka$ is 
$\al(X,v)(Y)=\dot b(Y)v-\si(X,Y)$.
With the notations ${\dot b}\check\ :v\in V\mapsto
\dot b(\cdot)v\in\Lin(\g,V)$ and $\si\check\ :X\in\g\mapsto i_X\si\in\Lin(\g,V)$, 
we write 
\begin{equation}\label{chec}
\al={\dot b}\check\ \o p-\si\check\ \o q. 
\end{equation}

To express the terms of  \eqref{inter} we need the identities: 
$$
\dot b(Y)\cdot p^{\eq}=\ev_Y\o(\dot b\check\ \o p)^{\eq}\text{ and }
q^*i_{Y^r}\si^{\eq}=-\ev_Y\o (\si\check\ \o q)^{\eq},
$$ 
where $\ev_Y:\Lin(\g,V)\to V$ for $Y\in\g$ denotes the evaluation map.
Indeed,
\begin{align*}
\dot b(Y)\cdot p^{\eq}(\hat g\cdot\hat X)
&=\dot b(Y)(b_g\cdot p(\hat X))
=b_g\cdot\dot b(\Ad_{g^{-1}}Y)p(\hat X)
=b_g\cdot\dot b\check\ (p(\hat X))(\Ad_{g^{-1}}Y)\\
&=(g\cdot (\dot b\check\ \o p)(\hat X))(Y)
={\ev}_{Y}\o(\dot b\check\ \o p)^{\eq}(\hat g\cdot\hat X)
\end{align*}
and 
\begin{align*}
q^*i_{Y^r}\si^{\eq}&(\hat g\cdot\hat X)
=\si^{\eq}(Y\cdot g,g\cdot X)
=b_g\cdot\si(\Ad_{g^{-1}}Y,X)
=-b_g\cdot(\si\check\ (X))(\Ad_{g^{-1}}Y)\\
&=-(g\cdot (\si\check\ \o q(\hat X)))(Y)
=-\ev_Y\o (\si\check\ \o q)^{\eq}(\hat g\cdot\hat X).
\end{align*}
We continue the computation \eqref{inter}:
\begin{align*}
\flux^{\eq}(Y^r,\dot b(Y))=\ev_Y\o[-(\si\check\ \o q)^{\eq}+({\dot b}\check\ \o p)^{\eq}]
\stackrel{\eqref{chec}}{=}\ev_Y\o[\al^{\eq}]=\ev_Y\o [d\ka]=0.
\end{align*}

Knowing now that $(\la,b)$ is an equi-hamiltonian action, we pull-back the prequantization abelian extension
\begin{equation*}
1\to A\stackrel{}{\to}\Diff^{\eq}(\hat G,\th=p^{\eq})^A_0
\stackrel{}{\to}\Diff_{\ex}^{\eq}(G,\om=-\si^{\eq})\to 1
\end{equation*}
by the equi-hamiltonian action $(\la,b)$. The result is an abelian Lie group extension of $G$ by $A$,
integrating the Lie algebra 2--cocycle $\si$ because
$$
-\om(\dot\la(X),\dot\la(Y))(e)=\si^{\eq}(X^r,Y^r)(e)=\si(X,Y),\quad X,Y\in\g.
$$
It is isomorphic with the given extension $\hat G$ since
they have the same Lie algebras and $G$ is simply connected and regular \cite{N04} Chapter VIII.
\end{proof}


\section{Group 2--cocycle}\label{s6}

In the case of an exact 2-form $\om=d\al$ for some $\al\in\Om^1(M, V)$, the principal bundle with curvature $\om$ 
is the trivial bundle $P=M\x V\stackrel{q}{\to}M$ 
with principal connection 1--form $\th=q^*\al+\th_V$, where $\th_V$ denotes the Maurer-Cartan form on $ V$.
Then the prequantization central extension (\ref{2222}) for $(M,d\al)$
is defined by $V$--valued group 2--cocycles on the hamiltonian group $\Diff_{\ex}(M,\om)$:
$$
c(\ph_1,\ph_2)=f(\ph_2)(x)-f(\ph_1\ph_2)(x)+f(\ph_1)(\ph_2(x))
$$
where $f:\Diff_{\ex}(M,\om)\to C^\oo(M,V)$ is a map satisfying $f(1_M)=0$ and $d(f(\ph))=\al-\ph^*\al$. Here $x\in M$ is arbitrary: $c(\ph_1,\ph_2)$ does not depend on $x$. 
The cocycle $c$ is cohomologous to the following group cocycle on $\Diff_{\ex}(M,\om)$:
$$
B(\ph_1,\ph_2)=\int_x^{\ph_2(x)}(\al-\ph_1^*\al).
$$
In the special case when $H^1(M,\RR)=0$, the cocycle $B$ can be extended to the group
$\Diff(M,\om)$ of $\om$-preserving diffeomorphisms  \cite{ILM06}.
The same holds for the cocycle $c$. This follows from the fact that, for any $\ph\in\Diff(M,\om)$,
the 1-form $\ph^*\al-\al$ is closed, hence exact.

We show that the prequantization abelian extension (\ref{four}) 
is described by a similar $V$-valued group 2--cocycle on $\Diff_{\ex}^{\eq}(M,\om)$.

\begin{theorem}\label{cB}
Given $\om=d\al$ for $\al\in\Om^1(M,V)$, 
the identity component of the group of equi-quantomorphisms $\Diff^{\eq}(P,\th)^A_0$
is the abelian extension of the group of equi-hamiltonian diffeomorphisms 
$\Diff_{\ex}^{\eq}(M,\om)$ by $V$ with cohomology class defined by the $V$--valued group 2--cocycle $c$ on $\Diff_{\ex}^{\eq}(M,\om)$,
\begin{equation*}\label{c}
c((\ph_1, u_1),(\ph_2, u_2))
= u_1(f(\ph_2, u_2)(x))-f(\ph_1\ph_2, u_1 u_2)(x)+f(\ph_1, u_1)(\ph_2(x)),
\end{equation*}
where the map
$f:\Diff_{\ex}^{\eq}(M,\om)\to C^\oo(M, V)$ is determined by
$f(1_M,1_ V)=0$ and $d(f(\ph, u))= u\cdot\al-\ph^*\al$.
Different choices for $f$ define cohomologous cocycles. 
The cocycle $c$ is cohomologous to the group cocycle
$$
B((\ph_1,u_1),(\ph_2,u_2))=\int_x^{\ph_2(x)}(u_1\cdot\al-\ph_1^*\al).
$$
\end{theorem}

\begin{proof}
The existence of the map $f$ follows from Remark \ref{deal}.

An arbitrary equi-quantomorphism  $\ps$ of $(P=M\x V,\th=q^*\al+\th_V)$
is of the form $\ps(x,v)=(\ph(x),m_\ps(x,v))$, where $\ph=q_*\ps$ and
$m_\ps:M\x V\to V$. The condition $\ps^*\th= u\cdot\th$ becomes
$dm_\ps=d(f(\ph, u)+ u)$, since $\ps^*(\th_V)=dm_\ps$ and $\ps^*q^*\al
=q^*\ph^*\al=q^*( u\cdot\al-d(f(\ph, u)))$. Hence the map $m_\ps$ is of the form $m_\ps(x,v)=f(\ph, u)(x)
+ u(v)+a$ for some $a\in V$ so
$$
\ps(x,v)=(\ph(x),f(\ph,u)(x)+u(v)+a).
$$ 
In this way we get a bijection $\ps\mapsto ((\ph, u),a)$
between $\Diff^{\eq}(P,\th)^A$ and the cartesian product
$\Diff^{\eq}_{\ex}(M,\om)\x V\subset\Diff(M)\x\GL(V)\x V$. 
The following computation shows 
that the group $\Diff^{\eq}(P,\th)_0^A$ is isomorphic to the
abelian extension defined by the given $V$--valued
group 2--cocycle $c$ on $\Diff_{\ex}^{\eq}(M,\om)$:
\begin{align*}
&(\ph_1, u_1,a_1)\o(\ph_2, u_2,a_2)(x,v)=
(\ph_1, u_1,a_1)(\ph_2(x),f(\ph_2, u_2)(x)+ u_2(v)+a_2)\\
&=(\ph_1\ph_2(x),f(\ph_1, u_1)(\ph_2(x))+ u_1f(\ph_2, u_2)(x)
+ u_1 u_2(v)+ u_1(a_2)+a_1)\\
&=(\ph_1\ph_2(x),f(\ph_1\ph_2, u_1 u_2)(x)+ u_1 u_2(v)+a_1+ u_1(a_2)
+c((\ph_1, u_1),(\ph_2, u_2)))\\
&=(\ph_1\ph_2,u_1u_2,a_1+ u_1(a_2)+c((\ph_1, u_1),(\ph_2, u_2)))(x,v)
\end{align*}
for all $(x,v)\in P$, which means 
$((\ph_1, u_1),a_1)\o((\ph_2, u_2),a_2)=((\ph_1,u_1)(\ph_2,u_2),a_1+ (\ph_1,u_1)\cdot a_2+c((\ph_1, u_1),(\ph_2, u_2)))$.

The cocycles $c$ and $B$ differ by the coboundary of the 1-cochain $(\ph,u)\mapsto f(\ph,u)(x)$ 
on $\Diff_{\ex}^{\eq}(M,\om)$,
so they are cohomologous and both describe the prequantization abelian extension.
\end{proof}

\begin{remark}
In the special case when $H^1(M,V)=0$, both cocycles $B$ and $c$ from Theorem \ref{cB} 
can be extended to the group $\Diff^{\eq}(M,\om)$ of $\om$-equivariant diffeomorphisms.
This follows from the fact that, for any $\ph\in\Diff^{\eq}(M,\om)$,
the $V$-valued 1-form $\ph^*\al-u\cdot\al$ on $M$ is closed, hence exact.
\end{remark}

\begin{remark}
Given an exact 2--form $\om=d\al\in\Om^2(M,V)$ and an equi-hamiltonian $G$--action $(\la,b)$, Theorem \ref{th2} provides an  abelian
Lie group extension $\hat G$ of $G$ by the $G$--module $V$, which can be defined also by a group 2--cocycle,
the pull-back of the cocycle from Theorem \ref{c} by the group homomorphism $(\la,b):G\to\Diff_{\ex}^{\eq}(M,\om)$:
$$
c(g_1,g_2)=b_{g_1}(f(g_2)(x))-f(g_1g_2)(x)+f(g_1)(\la_{g_2}(x)),
$$
where $f(g)\in C^\oo(M,V)$ with $f(e)=0$ and $d(f(g))=b_g\cdot\al-\la_g^*\al$, and $x\in M$ arbitrary.
The $\hat G$--action on $P=M\x V$ lifting the $G$--action and 
such that $\th$ is $\hat G$--equivariant
is $(g,v)\cdot (x,v')=(\la_g(x),f(g)(x)+b_g(v')+v)$.
\end{remark}


\section{Examples}\label{s7}


Theorem \ref{th2} provides a geometric construction of several 
abelian Lie group extensions of diffeomorphism groups.
In the examples below, the manifold $M$ will always be a homogeneous 
manifold $G/H$, with $H$ a connected Lie subgroup of $G$ and $V$ a $G$-module.
The $G$--equivariant closed 2--form $\om$ on $G/H$ is uniquely defined
by a $V$-valued Lie algebra 2--cocycle $\si$ on $\g$ satisfying two properties:
\begin{enumerate}
\item The kernel of $\si$ contains the Lie algebra $\h$ of $H$, so that $\si$ descends to a skew-symmetric bilinear form on $\g/\h$.
\item $\si$ is $H$--equivariant, \ie 
$\si(\Ad(g)X,\Ad(g)Y)=b(g)\cdot\si(X,Y)$ for all
$X,Y\in\g$ and $g\in H$.
The subgroup $H$ being connected, the last condition is equivalent
to 
$\si([Z,X],Y)+\si(X,[Z,Y])=\dot b(Z)\cdot\si(X,Y)$
for all $X,Y\in\g$ and $Z\in\h$.
\end{enumerate}

\begin{example}
We consider the group $\Diff_+(S^1)$ of orientation preserving 
diffeomorphisms of the circle and its right modules $\F_\la$ of
$\la$-densities on the circle: $b_\la(\ph)f=(\ph')^\la(f\o\ph)$
for $\ph\in\Diff_+(S^1)$ and $f\in C^\oo(S^1)$. 
The $\X(S^1)$--module structure on $\F_\la$ is given by
$\dot b_\la(X)f=Xf'+\la X'f$ for $X\in\X(S^1)$.
Here we identify $\la$-densities $f(x)(dx)^\la$ and vector fields $X(x)\frac{d}{dx}$ 
with smooth functions $f$ and $X$ on the circle.

The abelian extensions of $\X(S^1)$
defined with the $\F_\la$--valued cocycles 
\begin{gather*}
\si_0(X,Y)=\int_0^1(X'Y''-X''Y')dx\in\RR\subset \F_0\\
\si_1(X,Y)=X'Y''-X''Y'\in \F_1\\
\si_2(X,Y)=X'Y'''-X'''Y'\in \F_2
\end{gather*}
integrate to abelian extensions of $\Diff_+(S^1)$. Corresponding group
cocycles are presented in \cite{OR98}.
These abelian extensions can also be obtained geometrically 
by the Theorem \ref{th2}, 
taking $M$ to be the contractible homogeneous space $\Diff_+(S^1)/S^1$,
where $S^1$ is identified with the subgroup of rotations of $S^1$.
The existence of the $\Diff_+(S^1)$--equivariant 
$\F_\la$--valued closed 2--form $\om_\la$ on $M$ defined by the 2--cocycle $\si_\la$ for $\la=0,1,2$
is ensured by the $S^1$--equivariance of $\si_\la$
and the fact that the constant vector fields belong to the kernel of $\si_\la$.

The abelian extensions of $\X(S^1)$ by $\F_\la$ defined with the 2--cocycles  
\begin{gather*}
\bar\si_0(X,Y)=XY'-X'Y\in \F_0\\
\bar\si_1(X,Y)=XY''-X''Y\in \F_1\\
\bar\si_2(X,Y)=XY'''-X'''Y\in \F_2
\end{gather*}
integrate to abelian extensions of the universal covering group
$\widetilde\Diff_+(S^1)$ \cite{N04} Section 10. For the
geometric construction of these abelian extensions, in Theorem \ref{th2}
we take $G=\widetilde\Diff_+(S^1)$ acting by left
translations on $M=\widetilde\Diff_+(S^1)$, which is contractible.
The $\widetilde\Diff_+(S^1)$--equivariant $\F_\la$--valued
closed 2--form $\bar\om_\la$ is uniquely defined by its value
$\bar\si_\la$ at the identity, for $\la=0,1,2$.
\end{example}

\begin{example} 
Given a volume form $\mu$ on $M$, a non-trivial 
$\Om^1(M)/d\Om^0(M)$--valued Lie algebra 2--cocycle on
$\X(M)$ is $\si(X,Y)=(\div X)d(\div Y)$.
A group 2--cocycle on $\Diff(M)$ integrating $\si$ is constructed in \cite{B03}.

For the geometric construction
of an abelian Lie group extension of $\Diff(M)$ by 
its module $\Om^1(M)/d\Om^0(M)$, we remark that the 2--cocycle $\si$ is 
$\Diff(M,\mu)$--equivariant and the Lie algebra $\X(M,\mu)$ 
of divergence free vector fields is contained in the kernel of $\si$,
hence there is a closed $\Diff(M)$--equivariant 2--form $\om$ on the homogeneous 
space $\Diff(M)/\Diff(M,\mu)$. By a result of Moser \cite{KM97}, 
this space can be identified with the
contractible space of all volume forms of total mass 1.
Now Theorem \ref{th2} can be applied to $\om$.
\end{example}

\begin{example}
Let $\th$ be a connection 1--form on the principal $\GL(n,\RR)$-bundle
of frames $\pi:P(M)\to M$. Gelfand's cocycle presented in \cite{R06} is the 
$\Om^2(M)$--valued 2--cocycle $\si$ on $\X(M)$
defined by $\pi^*\si(X,Y)=\tr(L_{\tilde X}\th\wedge L_{\tilde Y}\th)$,
where $\tilde X,\tilde Y\in\X(P(M))$ are canonical lifts of $X,Y\in\X(M)$.
In the special case $M=\TT^n$ the $n$-torus, a group 2--cocycle on $\Diff(\TT^n)$ integrating $\si$
is constructed in \cite{B03}. 

For the 2--sphere $S^2$, Theorem \ref{th2} provides a geometric construction
of an abelian Lie group extension of $\Diff(S^2)$ by $\Om^2(S^2)$ integrating
Gelfand's cocycle.
We choose the connection $\th$ on the principal bundle of frames of $S^2$ 
induced by the canonical Riemannian metric on $S^2$.
The connected component of the isometry group of $S^2$ is $\SO(3)\subset\Diff(S^2)$,
so $L_{\tilde X}\th=0$ for $X\in\so(3)\subset\X(S^2)$.
The 2--cocycle $\si$ is $\SO(3)$--equivariant 
and the Lie algebra $\so(3)$  is contained in the kernel of $\si$,
hence there exists a closed $\Diff(S^2)$--equivariant 2--form $\om$ on the homogeneous 
space $\Diff(S^2)/\SO(3)$, given at the identity by $\si$. 
By a result of Smale \cite{S59}, this homogeneous space is contractible,
and Theorem \ref{th2} can be applied to $\om$.
\end{example}


\section{Appendix}

\subsection{Logarithmic derivative}

In this paragraph we collect some properties of the logarithmic derivative \cite{KM97}.

Given a manifold $M$ and a Lie group $G$ with Lie algebra $\g$,
the right logarithmic derivative of a function $h\in C^\oo(M,G)$
is $\de^r h\in\Om^1(M,\g)$:
$$
(\de^r h)( X_x):=(T_x h. X_x)h(x)^{-1},\forall X_x\in T_x M.
$$

\begin{remark}\label{rela}
A left logarithmic derivative $\de^l$ is defined by similar formula $(\de^l h)( X_x):=h(x)^{-1}(T_x h. X_x)$.
There is a relation between left and right logarithmic derivatives
\begin{equation}\label{relat}
\de^r h=\Ad(h)\de^l h=-\de^l(h^{-1}).
\end{equation}
\end{remark}

The right logarithmic derivative of a curve $h$ in the multiplicative group of positive real numbers is
$\de^r h(\partial_t)=\de^l h(\partial_t)=(\log h)'(t)$,
the derivative of the logarithm of $h$.

Given a curve $h$ in $G$ we will identify the $\g$-valued 1--form $\de^rh$ on $\RR$
with the curve $i_{\partial_t}\de^r h$ in $\g$.
When $\ph_t\in\Diff(M)$ is a diffeotopy of the manifold $M$,
then $\de^r\ph_t$ is the associated time dependent
vector field on $M$, and $\de^l\ph_t=\ph_t^*\de^r\ph_t$.
For $\om$ a differential form, we have $\tfrac{d}{dt}\ph_t^*\om=\ph_t^*L_{\de^r\ph_t}\om$.
This is a particular case of the identity
$
\tfrac{d}{dt}(g_t\cdot v)=\de^rg_t\cdot(g_t\cdot v)=g_t\cdot(\de^lg_t\cdot v),
$
for a smooth $G$--module $V$ and a curve $g_t$ in $G$.

\begin{remark}\label{rele}
Assume that $V$ is a Banach space.
Given a closed form $\om\in\Om^p(M,V)$, a path $\ph_t$ of diffeomorphisms of $M$ and a path $u_t$ of linear isomorphisms of $V$, both starting at the identity, the following equivalences hold for every $t$:
\begin{equation*}
\ph_t^*\om=u_t\cdot\om\Leftrightarrow L_{\de^r\ph_t}\om=\de^r u_t\cdot\om\Leftrightarrow L_{\de^l\ph_t}\om=\de^l u_t\cdot\om.
\end{equation*}
We show that $L_{\de^r \ph_t}\om=\de^r u_t\cdot\om$ implies
$\ph_t^*\om=u_t\cdot\om$.
The computation
$\frac{d}{dt}(\ph_t^*\om- u_t\cdot\om)=\ph_t^*L_{\de^r\ph_t}\om
-\de^r u_t\cdot (u_t\cdot\om)=\de^r u_t\cdot(\ph_t^*\om- u_t\cdot\om)$
shows that the curve
$\om_t=\ph_t^*\om- u_t\cdot\om$ is a solution of the differential equation
$\frac{d}{dt}\om_t=\de^r u_t\cdot\om_t$ in $\Om^p(M,V)$ with initial condition
$\om_0=0$. Evaluating at $p$ tangent vectors, we get a differential equation on the Banach space $V$,
for which uniqueness of solutions is ensured. It follows that $\om_t=0$, 
hence $\ph_t^*\om=u_t\cdot\om$.
\end{remark}

\begin{remark}\label{prod}
The right logarithmic derivative satisfies the Leibniz rule
\begin{equation*}
\de^r(h_1h_2)=\de^rh_1+\Ad(h_1)\de^rh_2
\end{equation*}
for $h_1,h_2\in C^\oo(M,G)$.
In other words $\de^r$ is a group 1--cocycle on $C^\oo(M,G)$
with values in the $C^\oo(M,G)$-module $\Om^1(M,\g)$ (for the adjoint action).
The logarithmic derivative of a $G$-valued function $h$ on a connected
manifold $M$ vanishes if and only if $h$ is constant. 
The Leibniz rule for the left logarithmic derivative is
\begin{equation}\label{dela}
\de^l(h_1h_2)=\de^l h_2+\Ad(h_2^{-1})\de^l h_1.
\end{equation}
\end{remark}

The {\it left Maurer-Cartan form} on a Lie group $G$ is the 1--form
$$
\th^l_G=\de^l(\id_G)\in\Om^1(G,\g).
$$ 
It satisfies $\th^l_G( X_g)=g^{-1}\cdot X_g$ for $X_g\in T_gG$,
and $\de^lh=h^*\th^l_G$ for any $h\in C^\oo(M,G)$.
Both $\th^l_G$ and $\de^lh$ satisfy the right Maurer-Cartan equation
$$
d\th^l_G+\frac12[\th^l_G,\th^l_G]=0,\quad
d\de^lh+\frac12[\de^lh,\de^lh]=0.
$$
Given $\al\in\Om^1(M,\g)$ which satisfies
the right Maurer-Cartan equation $d\al+\frac12\al\wedge\al=0$
and $U\subseteq M$ simply connected,
there exists a smooth function $h:U\to G$ with $\de^lh=\al$ on $U$.
In the special case of a domain in $\RR^2$,
when $h(t,s)\in G$ is a smooth two parameter family, we get 
\begin{equation}\label{ts}
\frac{d}{dt}\et-\frac{d}{dt}\xi=[\xi,\et],
\end{equation}
where $\xi(t,s)=(\frac{d}{dt}h(t,s))h(t,s)^{-1}$ and
$\et(t,s)=(\frac{d}{ds}h(t,s))h(t,s)^{-1}$.

\begin{proposition}\label{appendix}
Let $\th\in\Om^1(P,V)$ be a principal connection of the principal $A$-bundle $q:P\to M$ with $A=V/\Ga$
and principal action $\rho$. Then for any $f\in C^\oo(P,A)$,
\begin{equation}
\rh(f)^*\th=\th+\de^l(f),
\end{equation}
where $\rh$ denotes the principal $A$-action.
\end{proposition}

\begin{proof}
Let $q^{-1}(U)\cong U\x A$ be a principal bundle chart with $q$ and $q_A$ the two projections on $U$ and $A$. The principal connection 1--form in this chart is
$$
\th=q^*\al+q_A^*\th^l_A, 
$$
where $\th^l_A=\de^l(1_A)\in\Om^1(A,V)$ is the Maurer-Cartan form on $A$, and $\al\in\Om^1(U,V)$ closed.
The diffeomorphism $\rh(f)$ written in this chart is
$\rh(f)=(q,m\o(q_A,f))$, with $m$ denoting the multiplication map in $A$. Now
$$
\rh(f)^*\th=q^*\al+(m\o(q_A,f))^*\de^l(1_A)=q^*\al+\de^l(q_A)+\de^l(f)
=\th+\de^l(f),
$$
because by \eqref{dela} $\de^l(m\o(h_1,h_2))=\de^lh_1+\de^lh_2$ for $h_1,h_2:P\to A$.
\end{proof}


\subsection{1--Cocycles}\label{cite}

We list in this subsection some properties of 1--cocycles on Lie algebras and Lie groups \cite{N04}.
Let $G$ be a Lie group with Lie algebra $\g$ and let $V$ be a smooth $G$-module. Then $V$ is a $\g$-module and the pull-back action by the universal covering homomorphism $\tilde G\to G$ makes $V$ a $\tilde G$-module.

A $V$-valued {\it group 1--cocycle} is a locally smooth map $a:G\to V$ with
\begin{equation}\label{age}
a(gg')=a(g)+g\cdot a(g'),
\end{equation}
and a $V$-valued {\it Lie algebra 1--cocycle} is a continuous linear map $\al:\g\to V$ with
$$
\al([X,X'])=X\cdot\al(X')-X'\cdot\al(X).
$$
There is a natural map $a\mapsto \al=d_ea$ from locally smooth group 1--cocycles on the Lie group $G$ to Lie algebra 1--cocycles on its Lie algebra $\g$.

Let $\al^{\eq}$ be the equivariant $V$-valued 1-form on $G$ uniquely determined by the 1-cocycle $\al$ through $\al^{\eq}_g(g.X)=g\cdot \al(X)$ for all $g\in G$ and $X\in\g$. It is closed since $\al$ is a cocycle.

\begin{proposition}\label{properties}
If $a:G\to V$ is a group 1--cocycle integrating the Lie algebra 1--cocycle $\al:\g\to V$, then the following identities hold:
\begin{enumerate}
\item $da=\al^{\eq}\in\Om^1(G,V)$.
\item $\tfrac{d}{dt}(a(g_t))=g_t\cdot\al(\de^lg_t) =\al(\de^r g_t)+\de^r g_t\cdot a(g_t)$, for $g_t$ a path in $G$.
\end{enumerate}
\end{proposition}

\begin{proposition}\label{1cocy}
There exists a unique group 1--cocycle on the universal covering group $\tilde G$ integrating the Lie algebra 1--cocycle $\al:\g\to V$ 
$$
\tilde a:\tilde G\to V,\quad
\tilde a([g])=\int_0^1 g_t\cdot\al(\de^l g_t)dt,
$$ 
where $\tilde G$ is identified with the group of homotopy classes of piecewise smooth paths $g_t$ in $G$ starting at the identity.
\end{proposition}

\begin{proof}
Because $d\al^{\eq}=0$, the map $\tilde a$ is well defined.
The 1--cocycle condition \eqref{age} for $\tilde a$ is easily verified noticing that
the smooth path $g_tg'_t$, $t\in [0,1]$ and the piecewise smooth path $h_t$, defined by
$h_t=g_{2t}$ if $t\le\frac12$ and 
$h_t=g_1g'_{2t-1}$ if $t\ge\frac12$, are homotopic. 
\end{proof}

\begin{remark}
Let $\al:\g\to V$ be a Lie algebra 1-cocycle. 
If the group $\Pi_\al$ of periods of the closed 1--form $\al^{\eq}$ is discrete, 
then for any discrete subgroup $\Pi$ of $V$ containing $\Pi_\al$, 
the 1--cocycle $\tilde a$ descends to a 1--cocycle $a:G\to V/\Pi$. 
\end{remark}

From Proposition \ref{properties} and Remark \ref{rela} we deduce the following:

\begin{corollary}\label{kern}
Let $a:G\to V/\Pi$ be a group 1--cocycle integrating the Lie algebra 1--cocycle $\al:\g\to V$.
Then $\Ker a$ is a subgroup of $G$ and $\Ker \al$ is a Lie subalgebra of $\g$ such that, for any smooth curve $g$ in $G$ starting at the identity, the following are equivalent:
$g_t\in\Ker a$$\Leftrightarrow$$\de^l g_t\in\Ker\al$$\Leftrightarrow$$\de^r g_t\in\Ker\al$.
\end{corollary}

\begin{remark}\label{refinement}
The results of this section about 1-cocycles hold also for diffeomorphism groups $G$ with Lie algebras 
$\g$ of vector fields such that  
the following are equivalent:
$g_t\in G$$\Leftrightarrow$$\de^l g_t\in\g$$\Leftrightarrow$$\de^r g_t\in\gl$.
\end{remark}


\subsection{Obstructions to integrability}

According to the general theory developed in \cite{N02} and 
\cite{N04}, there 
are two obstructions for the integration of a Lie algebra 
cocycle $\si$ on $\g$ with values in a 
$G$-module $V$ to a Lie group extension of $G$ by a quotient group of $V$:  
the period map and the flux homomorphism. 

Let $\si^{\eq}$ be the closed equivariant $V$-valued 2-form on $G$ uniquely determined 
by the 2-cocycle $\si$. 
The {\it period map} is the group homomorphism 
$$\per_\si:\pi_2(G)\to V^G, \quad 
\per_\si([ c])=\int_{S^2} c^*\si^{\rm eq} \quad 
\mbox{ for } \quad  c\in C^\infty(S^2,G).$$ 
Its image $\Pi_\si$ 
is called the {\it period group} of $\si$. 

The {\it flux homomorphism} $F_\si:\pi_1(G)\to H^1(\g,V), [\gamma] 
\mapsto [I_\gamma^\si]$, 
assigns to each piecewise smooth 
loop $\ga$ in $G$ based at the identity,
the cohomology class of the 1--cocycle
$$I_\ga^\si : \g\to V, \quad I_\gamma^\si(X) = -\int_\ga i_{X^r}\si^{\rm eq}.$$

\begin{theorem} 
\label{integrate}
For a Lie algebra $V$-valued 2--cocycle $\si$ on $\g$  with discrete
period group $\Pi_\si$ and vanishing flux homomorphism $F_\si$,
the abelian Lie algebra extension $0\to V\to\hat\g\to\g\to 0$ defined by $\si$
integrates to an abelian Lie group extension
$$
1\to V/\Pi_\si\to\hat G\to G\to 1.
$$ 
\end{theorem}



\paragraph{Acknowledgements.}
The author thanks Karl-Hermann Neeb for many very useful comments and suggestions.
The valuable remarks and corrections due to the anonymous referee, especially those concerning 
Theorem \ref{th1}, are gratefully acknowledged. 



\end{document}